\input amstex
\documentstyle{amsppt}
\magnification=\magstep1

\NoBlackBoxes
\TagsAsMath

\pageheight{9.0truein}
\pagewidth{6.5truein}

\long\def\ignore#1\endignore{\par DIAGRAM\par}
\long\def\ignore#1\endignore{#1}

\ignore
\input xy \xyoption{matrix} \xyoption{arrow}
          \xyoption{curve}  \xyoption{frame}
\def\edge{\ar@{-}}
\def\dttdar{\ar@{.>}}

\def\dashedge{\ar@{--}}

\def\dshdar{\ar@{-->}}
\def\dropdown#1{\save+<0ex,-4ex> \drop{#1} \restore}

\def\dbull{*{\bullet}+0}
\def\mbull{\ar@{}[r]|{\bullet}}

\def\loopNE{\ar@'{@+{[0,0]+(6,2)} @+{[0,0]+(10,10)}
@+{[0,0]+(2,6)}}}
\def\loopNW{\ar@'{@+{[0,0]+(-2,6)} @+{[0,0]+(-10,10)} 
@+{[0,0]+(-6,2)}}}
\def\loopSW{\ar@'{@+{[0,0]+(-6,-2)} @+{[0,0]+(-10,-10)}
@+{[0,0]+(-2,-6)}}}
\def\loopSE{\ar@'{@+{[0,0]+(2,-6)} @+{[0,0]+(10,-10)}
@+{[0,0]+(6,-2)}}}

\def\loopNNE{\ar@'{@+{[0,0]+(4,2)} @+{[0,0]+(6,11)}
@+{[0,0]+(0,6)}}}
\def\loopSSW{\ar@'{@+{[0,0]+(-4,-2)} @+{[0,0]+(-6,-12)}
@+{[0,0]+(0,-6)}}}
\def\loopSSE{\ar@'{@+{[0,0]+(0,-6)} @+{[0,0]+(6,-11)}
@+{[0,0]+(4,-2)}}}
\endignore

\def\vsubseteq{\hbox{$\bigcup$\kern0.1em\raise0.05ex\hbox{$\tsize|$}}}
\def\smvsubseteq{\hbox{$\ssize\bigcup$\kern0.03em\raise0.05ex\hbox{$\ssize|$}}}

\def\la{{\Lambda}}
\def\lamod{\Lambda\text{-}\roman{mod}}
\def\modla{\roman{mod}\text{-} \Lambda}

\def \len{\operatorname{length}} 
\def\AA{{\Bbb A}}
\def\CC{{\Bbb C}}

\def\SS{{\Bbb S}}

\def\NN{{\Bbb N}}

\def\hom{\operatorname{Hom}}
\def\aut{\operatorname{Aut}}
\def\Aut{\operatorname{Aut}}

\def\soc{\operatorname{soc}}
\def\radsocbd{\operatorname{\bold{rad-soc(\bd)}}}

\def\ann{\operatorname{ann}}
\def\nullity{\operatorname{nullity}}

\def\ann{\operatorname{ann}}

\def\GL{\operatorname{GL}}

\def\Im{\operatorname{Im}}

\def\Seq{\operatorname{\frak{Seq} }}
\def\rep{\operatorname{Rep}}

\def\Ext{\operatorname{Ext}}

\def\End{\operatorname{End}}

\def\udim{\operatorname{\underline{dim}}}

\def\A{{\Cal A}} 
 
\def\C{{\Cal C}}
\def\D{{\Cal D}}

\def\T{{\Cal T}}

\def\S{{\sigma}}

\def\m{{\frak m}}
\def\T{{\Cal T}}

\def\U{{\frak U}}
\def\V{{\frak V}}

\def\bd{\bold {d}}

\def\Qhat{\widehat{Q}}

\def\SStilde{\widetilde{\SS}}
\def\Btilde{\widetilde{B}}
\def\btilde{\widetilde{b}}
\def\Gtilde{\widetilde{G}}

\def\Ktilde{\widetilde{K}}

\def\latilde{\widetilde{\Lambda}}

\def\Mtilde{\widetilde{M}}

\def\bP{\bold{P}}

\def\Qhat{\widehat{Q}}

\def\bz{\bold{z}}

\def\modlad{\operatorname{\bold{Rep}}_{\bold{d}}(\Lambda)}
\def\modlasd{\operatorname{\bold{Rep}}_{d}(\Lambda)}

\def\Mod{\operatorname{\bold{Rep}}}

\def\toptd{\operatorname{\bold{Rep}^T_{\bold{d}}}}

\def\laySS{\operatorname{\bold{Rep}} \SS}

\def\laySSprime{\operatorname{\bold{Rep}} \SS'}

\def\Hom{\operatorname{Hom}}

\def\Gr{\operatorname{Gr}}

\def\grassS{\operatorname{\frak{Grass}}(\S)}

\def\grassSS{\operatorname{\frak{Grass}} \SS}

 \def\biggrass{\operatorname{GRASS}}
  \def\biggrasslasd{\operatorname{GRASS}_{d}(\la)}
 \def\biggrasslad{\operatorname{GRASS}_{\bd}(\la)}

\def\Schu{\operatorname{Schu}}

\def\Schu{\operatorname{\Schu}}

\def\GRASS{\operatorname{GRASS}}

\def\grassbd{\GRASS_{\bold d}(\Lambda)}

\def\Gr{\operatorname{Gr}}

\def\grassS{\operatorname{\frak{Grass}}(\S)}

\def\grassSS{\operatorname{\frak{Grass}} \SS}
\def\biggrassSS{\GRASS(\SS)}

\def\biggrassS{\GRASS(\S)}

\def\biggrass{\GRASS}

\def\Schu{\operatorname{Schu}}

\def\BHT{{\bf 1}}
\def\BCH{{\bf 2}}
\def\BoHZtwo{{\bf 3}}
\def\CW{{\bf 4}}
\def\CBS{{\bf 5}}
\def\DoFl{{\bf 6}}
\def\EiSa{{\bf 7}}
\def\Ger{{\bf 8}}
\def\Gur{{\bf 9}}
\def\hier{{\bf 10}}
\def\GoHu{{\bf 11}}
\def\KacI{{\bf 12}}
\def\KacII{{\bf 13}}
\def\Kra{{\bf 14}}
\def\Mor{{\bf 15}}
\def\RiRuSm{{\bf 16}}
\def\Scho{{\bf 17}}
\def\Schro{{\bf 18}}

\topmatter

\hphantom{}
\vskip-0.7truein

\title Irreducible components of varieties of representations I.  The local case. \endtitle

\rightheadtext{irreducible components}

\author B. Huisgen-Zimmermann
\endauthor

\thanks The author was partially supported by an NSF grant while carrying out this work.  Moreover, she was supported by NSF grant 0932078 000, while in residence at MSRI, Berkeley. \endthanks

\address Department of Mathematics, University of California, Santa
Barbara, CA 93106-3080 \endaddress
\email birge\@math.ucsb.edu \endemail

\abstract  Let $\la$ be a local truncated path algebra over an algebraically closed field $K$, i.e., a quotient $KQ/ \langle \text{the paths of length \, } L+1 \rangle$ of a path algebra $KQ$, where $L$ is a positive integer and $Q$ is the quiver with a single vertex and a finite number $r$ of loops.  For any positive integer $d$, we determine the irreducible components of the varieties that parametrize the $d$-dimensional representations of $\la$, namely, the components of the classical affine variety $\modlasd$ and  --  equivalently  --  those of the projective parametrizing variety $\biggrasslasd$.  The former variety consists of  the $r$-tuples $(A_1, \dots, A_r)$ of $d \times d$-matrices subject to the condition that all $(L+1)$-fold products of the entries vanish, while $\biggrasslasd$ is the closed subvariety of the classical Grassmann variety $\Gr\bigl((\dim \la^d) - d, \, \la^d\bigr)$ which picks out the $\la$-sub{\it modules\/} of codimension $d$ of the free left $\la$-module $\la^d$.  Our method is to corner the components by way of a twin pair of upper semicontinuous maps from $\modlasd$ to a poset consisting of sequences of semisimple modules.

An excerpt of the main result is as follows.  Given a sequence $\SS = (\SS_0, \dots, \SS_L)$ of semisimple modules with $\dim \bigoplus_{0 \le l \le L} \SS_l = d$, let $\laySS$ be the subvariety of $\modlasd$ consisting of the points that parametrize the modules with radical layering $\SS$. (The radical layering of a $\la$-module $M$ is the sequence $\bigl(J^l M / J^{l+1} M\bigr)_{0 \le l \le L}$, where $J$ is the Jacobson radical of $\la$.)  Suppose the quiver $Q$ has $r \ge 2$ loops.  If $d \le L+1$, the variety $\modlasd$ is irreducible and, generically, its modules are uniserial.  If, on the other hand, $d > L+1$, then the irreducible components of $\modlasd$ are the closures of the subvarieties $\laySS$ for those sequences $\SS$ which satisfy the inequalities $\dim \SS_l \le r \dim \SS_{l+1}$ and $\dim \SS_{l+1} \le r \dim \SS_l$ for $0 \le l < L$; generically, the modules in any such component have socle layering $(\SS_L, \dots, \SS_0)$.  As a byproduct, the main result provides further installments of generic information on the modules corresponding to the irreducible components of the parametrizing varieties. 
\endabstract
   
\endtopmatter

\document

\head  Introduction and main result  \endhead 

The broader objective of this article is to further develop a strategy of generically organizing the representation theory of basic finite dimensional algebras $\la$ over an algebraically closed field $K$; those are the path algebras modulo relations, $\la = KQ/I$, where $Q$ is a quiver and $I$ an admissible ideal.  The strategy amounts to tackling the following two sequential tasks for each dimension vector $\bd$ of $\la$:  (I)  Determine, in both geometric and representation-theoretic terms, the irreducible components of the affine parametrizing varieties $\modlad$ for the $\bd$-dimensional $\la$-modules, and (II) explore generic features of the modules corresponding to the individual components, such as generic decomposition properties, generic tops, socles, submodule lattices, etc.  

This approach is not new.  It is rooted in work of Kac (see \cite{\KacI}, \cite{\KacII}) and Schofield (see \cite{\Scho}) focused on hereditary algebras, i.e., on the case where $I = 0$.  In this situation, task (I) is non-existent, in that the varieties $\modlad$ are always full affine spaces.  In attempts to extend the generic analysis of the modules with fixed dimension vector to non-hereditary algebras, the first hurdle encountered is the fact that the  $\modlad$ split into a plethora of irreducible components in general.   Hence, the goal of determining these components from $Q$ and $I$ was singled out early on, by Kraft in \cite{\Kra} at the latest.  Results of Crawley-Boevey and Schr\"oer targeted those irreducible components of $\modlad$ whose representations are generically decomposable, relating them to components for smaller dimension vectors (see \cite{\CBS}).  Moreover, Babson, Thomas and the author showed how to obtain (for arbitrary choices of $Q$ and $I$) a finite list of representation-theoretically defined closed irreducible subvarieties of $\modlad$ which includes the irreducible components (see \cite{\BHT, Sections 3-5}); this turns the component problem into a sifting problem.  For vanishing radical square, task (I) was recently completed by Bleher, Chinburg and the author, and task (II) was pushed to the level achieved for hereditary algebras (see \cite{\BCH}).   Beyond that, a full description of the irreducible components is available for certain special biserial algebras, notably for the family of Gelfand-Ponomarev algebras $K[X, Y]/ (XY, X^r, Y^s)$ (see \cite{\Schro}), next to the algebra $K[X,Y]/ (X^2, Y^2)$ (see \cite{\RiRuSm}) and special biserial algebras dubbed {\it gentle\/} (see {\cite{\CW}).  These algebras have tame representation type, and explicit descriptions of their finitely generated indecomposable modules help access the components of their parametrizing varieties. 

Two types of parametrizing varieties are under consideration, the classical affine varieties $\modlad$ and their projective counterparts $\grassbd$.  Tasks (I) and (II) are equivalent in the two scenarios, in that any information gained in one can readily be transferred to the other (see Proposition 1.2); but the two settings draw on different caches of geometric techniques, making it advantageous to switch back and forth.   The varieties $\modlad$ and $\grassbd$ are connected, and the abbreviated terminology ``components" will always refer to ``irreducible components".  When $\la$ is local, meaning that $Q$ has only one vertex, the dimension vector $\bd$ of a module is replaced by its dimension $d$.   

In light of the combinatorial intricacies surfacing in known examples, it does not appear reasonable to aim for a meaningful classification of the components in a ``closed format", that is, in a format covering arbitrary choices of $\la$.  Rather, the problem calls for an understanding of ``base cases" and an arsenal of techniques for departing from those.   The class of departure we propose consists of the {\it truncated path algebras\/}, i.e., the algebras
$\la = KQ/ \langle \text{all paths of length\ } L+1 \rangle$ for some $L \ge 1$.  These algebras  --  they include the hereditary ones and those with vanishing radical square  --  hold a prominent position among the basic finite dimensional $K$-algebras $\la'$.  Indeed, each $\la'$ is isomorphic to a factor algebra of a unique truncated path algebra $\la$ with the same quiver and Loewy length, whence each $\Mod_{\bd}(\la')$ arises as a closed subvariety of the corresponding $\modlad$. This makes the truncated case a natural first target in developing the suggested approach. 

In the truncated situation, the following subvarieties $\laySS$ are irreducible and the (finite) family of their closures includes all irreducible components of $\modlad$ (see \cite{\BHT}).  We recall the relevant definitions: The {\it radical layering\/} of a module $M$ is the sequence $\SS(M) = \bigl(J^l M/ J^{l+1}M\bigr)_{0 \le l \le L}$, and the {\it socle layering\/} $\SS^*(M)$ is defined dually. Both layerings are {\it semisimple sequences with dimension vector $\udim M$\/}, i.e., sequences of the form $\SS = (\SS_0, \SS_1, \dots, \SS_L)$ whose entries $\SS_l$ are semisimple modules such that $\udim \SS : = \sum_{0 \le l \le L} \udim \SS_l = \udim M$.  The variety $\laySS$ is the locally closed subset of $\modlad$ consisting of the points which parametrize the modules with radical layering $\SS$.   
Consequently, Task (I) above amounts to characterizing those semisimple sequences $\SS$ which arise as the generic radical layerings of the irreducible components of $\modlad$.  For the purpose of filtering out the relevant sequences $\SS$, we consider the ``dominance" partial order on the set of semisimple sequences (see under conventions below), and use the following observation as a starting point (see Theorem 3.1).   

\proclaim{Observation} Let $\la$ be truncated.  If $\, \SS = \SS(M)$ for some module $M$ with the property that the pair $\bigl(\SS(M), \SS^*(M) \bigr)$ is a minimal element of the set

\centerline{$ \radsocbd :=   \{ \bigl(\SS(M), \SS^*(M) \bigr) \mid M \in \lamod, \, \udim M = \bd \}$,}

\noindent then the closure $\overline{\laySS}$ is an irreducible component of $\modlad$.  
\endproclaim

In general, the converse fails, as witnessed by Example 4.5.  But in a number of significant special cases it does hold, notably when $\la$ is local, meaning that the quiver $Q$ has the form

$$\xymatrixrowsep{0.75pc}\xymatrixcolsep{0.3pc}
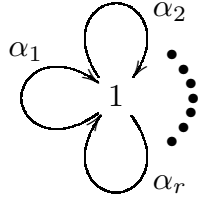
\xymatrix{
 & \ar@{{}{*}}@/^1pc/[dd] \\
1 \ar@'{@+{[0,0]+(-6,-6)} @+{[0,0]+(-15,0)} @+{[0,0]+(-6,6)}}^(0.7){\alpha_1}
\ar@'{@+{[0,0]+(-6,6)} @+{[0,0]+(0,15)} @+{[0,0]+(6,6)}}^(0.7){\alpha_2}
\ar@'{@+{[0,0]+(6,-6)} @+{[0,0]+(0,-15)} @+{[0,0]+(-6,-6)}}^(0.3){\alpha_r}  \\
 &
}$$

\noindent for some $r \in \NN$.  

In case $r = 1$,  $\la$ is a truncated polynomial ring in a single variable, and all of the varieties $\modlasd$ are trivially irreducible.  Otherwise, we find:

\proclaim{Theorem A}  
Let $\la$ be a local truncated path algebra of Loewy length $L+1$, whose quiver $Q$ has $r \ge 2$ loops, and let $d$ be a positive integer.  

If $d > L+1$, the irreducible components of $\modlasd$ are precisely the subvarieties $\overline{\laySS}$, where $\SS$ is a semisimple sequence of dimension $d$ satisfying the following equivalent conditions:
\roster
\item The closure $\overline{\laySS}$ is an irreducible component of $\modlasd$.
\smallskip
\item"(1')" The closure $\, \overline{\biggrass\, \SS}$ is  an irreducible component of $\biggrass_d(\la)$.
\smallskip
\item"(2)"   $\dim \SS_l \le r \cdot \dim \SS_{l-1}$ and $\ \dim \SS_{l-1} \le r \cdot \dim \SS_l\,$ for $\, l \in \{1, \dots, L\}$; in particular, $\SS_l  \ne 0$ for all $l \le L$.
\smallskip
\item"(3)" $\laySS \ne \varnothing$, and $\SS^* = (\SS_L, \SS_{L-1}, \dots, \SS_0)$ is the generic socle layering of the modules in $\laySS$.
\smallskip
\item"(4)"   $\SS = \SS(M)$ for some minimal pair $\bigl(\SS(M), \SS^*(M) \bigr)$  in $\radsocbd$.
\endroster
\smallskip
If, on the other hand, $d \le L+1$, the variety $\modlasd$ is irreducible and, generically, its modules are uniserial.
\endproclaim

We remark that, for $L = 1$ and $r = 2$, the irreducible components of the $\modlasd$ were already determined by Donald and Flanigan \cite{\DoFl}, as well as by Morrison \cite{\Mor}.  For arbitrary choice of $r$, the case $L=1$ is covered by \cite{\BCH, Theorem 3.12}.

Condition (2) of the theorem permits us to list the irreducible components of $\modlasd$, tagged by their generic radical layerings, for all values of $r$, $L$ and $d$.  Beyond the modest installment of progress towards Task (II) stated explicitly, the theorem provides some further generic mileage (Section 4, Supplement B and Corollaries C-E). 

In a subsequent paper, Babson, Shipman, and the author will show that the equivalence ``(1)$\iff$(4)" extends to the situations where $\la$ is truncated and either based on an acyclic quiver $Q$ or else has Loewy length at most $3$ (the latter bound being sharp).  We deal with these cases separately, since they are fraught with their own technical difficulties.  

Section 2 applies to arbitrary basic finite dimensional algebras.  It provides the framework for our technique of cornering the irreducible components of $\modlad$ and $\biggrasslad$ by way of upper semicontinuous maps.  Section 3 narrows the focus to truncated path algebras, preparing for both applications and a proof of Theorem A.  In Section 4, we derive consequences from the theorem and buttress the theory with examples.  A proof of Theorem A is given in Section 5. 

Finally, we point out that interest in varieties which consist of sequences of matrices satisfying certain relations is not limited to the role they play in the representation theory of algebras; see, e.g., \cite{\Ger}, \cite{\EiSa}, \cite{\Gur}.
\bigskip

\head 1. Conventions and prerequisites \endhead 

 Throughout, $\la$ is a basic finite dimensional algebra over an algebraically closed field $K$, and $J$ denotes its Jacobson radical; say $J^{L+1} = 0$.  Without loss of generality, we will assume that $\la  = KQ/ I$ for a quiver $Q$ and some admissible ideal $I$ in the path algebra $KQ$.  Moreover, the set $Q_0 = \{e_1, \dots, e_n\}$ of vertices of $Q$ will be identified with a full set of primitive idempotents of $\la$.  

An element $x \in M$ is said to be {\it normed\/} if $x = e_i x$ for some $i$.  A {\it top element\/} of $M$ is a normed element  in $M \setminus JM$, and a {\it full sequence of top elements of $M$\/} is any generating set of $M$ consisting of top elements which are $K$-linearly independent modulo $JM$. 
\smallskip

\noindent{\bf Semisimple sequences:}
By $S_i$ we denote the simple module $\la e_i/ J e_i$.  
Unless we want to distinguish among different embeddings, we identify isomorphic semisimple modules throughout.  This identification provides us with a partial order on the set of finite dimensional semisimple modules:  Namely $U \subseteq V$  if and only if the multiplicity of each $S_i$ in $U$ is smaller than or equal to the multiplicity of $S_i$ in $V$.
 
A {\it semisimple sequence\/} in $\lamod$ is a sequence of the form $\SS = (\SS_0, \SS_1, \dots, \SS_L)$ such that each entry $\SS_l$ is a semisimple module.  The {\it dimension vector of\/} $\SS$, written $\udim \SS$, is defined to be $\sum_{0 \le l \le L}\udim \SS_l$.   The set of all semisimple sequences with fixed dimension vector $\bd$ is endowed with the following partial order, dubbed the {\it dominance order\/}:  
$$\SS \le \SS' \ \ \ \iff \ \ \ \bigoplus_{0 \le j \le l} \SS_j \subseteq \bigoplus_{0 \le j \le l} \SS'_j \ \ \text{for}\  \ l \in \{0, 1, \dots, L\}.$$

Two types of semisimple sequences will play a pivotal role in the sequel, namely the radical and socle layerings of the finitely generated $\la$-modules:  For $M \in \lamod$, the {\it radical layering\/} of $M$ is $\SS(M) = (J^l M/J^{l+1}M)_{0 \le l \le L}$.  Dually, we define the {\it socle layering\/} of $M$ to be the sequence $\SS^*(M) = \bigl(\soc_l M/ \soc_{l-1} M \bigr)_{0 \le l \le L}$,  where $\soc_{-1} M = 0$ and $\soc M = \soc_0 M \subseteq \soc_1 M \subseteq \cdots \subseteq \soc_L M$ is the standard socle series of $M$:  Namely, $\soc_0 M = \soc M$, and 
$$\soc_{l+1} M / \soc_l M = \soc (M/\soc_l M).$$
Note that radical and socle layerings are dual to each other, in the sense that
$$\SS(D(M)) = \bigl(D(\SS^*_0(M)), \cdots, D(\SS^*_L(M) \bigr)  \  \text{and}\, \  \SS^*(D(M)) =  \bigl(D(\SS_0(M)), \cdots, D(\SS_L(M) \bigr);$$ 
here $D$ denotes the duality $\Hom_K( - ,K)$ between the categories $\lamod$ and $\modla$ of finitely generated left, resp\. right, $\la$-modules.
Clearly,  $\udim\, \SS(M) = \udim\, \SS^*(M) = \udim M$.

Our interest in semisimple sequences is restricted to those which arise as radical or socle layerings of $\la$-modules.   Whenever we have the choice, we will prioritize radical layerings, whence the bias in the following definition.

\definition{Definition 1.1}  A semisimple sequence $\SS$ is said to be {\it realizable\/} if there exists a left $\la$-module $M$ with $\SS(M) = \SS$.
\enddefinition

\noindent{\bf The parametrizing varieties:}  We recall the definitions of the relevant varieties parametri\-zing the isomorphism classes of (left) $\la$-modules with dimension vector $\bd = (d_1, \dots, d_n)$.    
\smallskip

\noindent  {\it The affine setting\/}:  The classical affine variety  is $\modlad = $
$$ \bigl\{ (x_\alpha)_{\alpha \in Q_1} \in \prod_{\alpha \in Q_1} \Hom_K \bigl(K^{d_{\text{start}(\alpha)}},\, K^{d_{\text{end}(\alpha)}}\bigr) \mid \text{the\ } x_{\alpha} \ \text{satisfy all relations in}\ I \bigr\},$$
where $Q_1$ is the set of arrows of the quiver $Q$.
This variety carries a conjugation action by $\GL(\bd) := \GL_{d_1}(K) \times \cdots\times \GL_{d_n}(K)$, the orbits of which are in bijective correspondence with the isomorphism classes of modules having dimension vector $\bd$.  Throughout, we denote by $M_x \in \lamod$ the module that corresponds to a point $x \in \modlad$.
If $\SS$ is a semisimple sequence with this dimension vector, $\laySS$ stands for the locally closed subvariety consisting of those points $x \in \modlad$ which have the property that the corresponding module $M_x$ has radical layering $\SS$.  Note that $\laySS \ne \varnothing$ if and only if $\SS$ is realizable in the sense of Definition 1.1.
\smallskip

\noindent  {\it The projective setting\/}: Let $d = |\bd| = \sum_i d_i$.  We fix a projective module $\bP = \bigoplus_{1 \le r \le d} \la \bz_r$ whose top has dimension vector $\udim\, (\bP / J \bP) = \bd$;  here $\bz_1, \dots, \bz_d$ is a full sequence of top elements of $\bP$.  In other words, $\bP$ is a projective cover of $\,\bigoplus_{1 \le i \le n} S_i^{d_i}$, and thus is the smallest projective $\la$-module with the property that every module with dimension vector $\bd$ is isomorphic to a quotient $\bP/C$ for some submodule $C \subseteq \bP$.  The variety $\Gr\bigl( (\dim \bP - d), \bP\bigr)$ is the vector space Grassmannian of all $(\dim\bP - d)$-dimensional $K$-subspaces of $\bP$, and $\biggrasslad$ is the closed subset consisting of the $\la$-submodules $C$ of $\bP$ with the property that $\udim \bP/C = \bd$.  Under the canonical action of the automorphism group $\Aut_\la(\bP)$ on $\biggrasslad$, the orbits are again in one-to-one correspondence with the isomorphism classes of modules with dimension vector $\bd$.  In parallel with the affine setting:  For any semisimple sequence $\SS$ with dimension vector $\bd$, we denote by $\biggrassSS$ the locally closed subvariety of $\biggrasslad$ picking up the points $C$ with $\SS(\bP/C) = \SS$.  (Caveat:  The variety $\biggrassSS$ introduced here is not to be confused with the much smaller one, $\grassSS$, used in \cite{\BHT}, for instance;  it is in this smaller variety that information on $\biggrassSS$ is preferably gleaned.)  
\smallskip

\noindent {\it Connection between the two settings\/}:
The horizontal double arrows in the diagram below point to the transfer of geometric information spelled out in the upcoming proposition.   It was proved in \cite{\BoHZtwo}, modulo the unirationality statement which was added in \cite{\GoHu}. 

\ignore
 \goodbreak \midinsert
$$\xymatrixrowsep{0.2pc}\xymatrixcolsep{3pc}
\xymatrix{
\boxed{\grassbd} \dropdown{\txt{(projective)}} \ar@{<->}[r] &\boxed{\modlad}
\dropdown{\txt{(affine)}} \\  \\
\smvsubseteq &\smvsubseteq \\
\boxed{\biggrassSS} \dropdown{\txt{(quasi-projective)}} \ar@{<->}[r]
&\boxed{\laySS} \dropdown{\txt{(quasi-affine)}}
 }$$
\medskip
\centerline{Diagram 1.1}
\endinsert 
\endignore

\proclaim{Proposition 1.2. Information transfer between the affine and projective settings}   
Consider the one-to-one correspondence between the orbits of $\, \grassbd$ and $\, \modlad $ assigning to any orbit $\aut_\la(\bP).C \subseteq \grassbd$ the orbit
$\GL(\bd).x \subseteq \modlad $ that represents the same isomorphism class of
$\la$-modules.  This correspondence extends to an inclusion-pre\-serv\-ing bijection
$$\Phi: \{ \aut_\la(\bP)\text{-stable subsets of\ } \grassbd \} \rightarrow
\{\GL(\bd)\text{-stable subsets of\ } \modlad \}$$  
which preserves and reflects
openness, closures, irreducibility, smoothness and unirationality. 

 In particular,  a semisimple sequence $\SS$ is the generic radical layering of the modules parametrized by an irreducible component of $\modlad$ precisely when $\SS$ has the same property relative to an irreducible component of $\, \grassbd$.
\endproclaim

For more detail regarding the parametrizing varieties, see \cite{\hier}.
\smallskip

\noindent{\bf Further terminology:} For any subset $\U$ of $\modlad$ or $\biggrasslad$, we refer to the modules corresponding to the points in $\U$ as the modules ``in" $\U$.  When $\U$ is irreducible, the modules in $\U$ are said to {\it generically\/} have property $(*)$ in case all modules in some dense open subset of $\U$ satisfy $(*)$.  Observe that, given any irreducible subset $\U \subseteq \modlad$, the radical layerings and socle layerings of the modules in $\U$ are generically constant (we re-emphasize our convention of identifying isomorphic semismple modules), whence it is meaningful to speak of {\it the\/} generic radical and socle layerings of the irreducible components of $\modlad$.

\head 2.  Strategy:  Cornering the components of $\modlad$ via upper semicontinuous maps \endhead

Throughout this section, we let $\la = KQ/I$ be an arbitrary basic finite dimensional algebra.   

\definition{Definition 2.1} Suppose $X$ is a topological space and $(\A, \le )$ a poset.  For $a \in \A$, we denote by $[a, \infty)$ the set $\{b \in \A \mid b \ge a\}$; the sets $(a, \infty)$, $(- \infty, a]$ and  $(- \infty , a)$ are defined analogously.

A map 
$f: X \longrightarrow \A$ is called {\it upper semicontinuous\/} if, for every element $a \in \A$, the pre-image of $[a, \infty)$ under $f$ is closed in $X$.  
\enddefinition

The following module invariants are well-known to yield upper semicontinuous maps on $X = \modlad$.  They all take numerical values and have finite images, hence satisfy the hypotheses of the upcoming observation. For any fixed $N \in \lamod$, the maps $x \mapsto \dim \hom_\la(M_x, N)$ and $x \mapsto \dim \hom_\la(N, M_x)$, $x \mapsto \dim \Ext^1_\la(M_x, N)$, and $x \mapsto \ \dim \Ext^1_\la(N,M_x)$ are examples; for $\Ext^1$, see \cite{\CBS}.  Moreover, for any path $p$ in $KQ \setminus  I$, the map $x \mapsto \nullity_p M_x$ is upper semicontinuous; here $\nullity_p M_x$ is the nullity of the $K$-linear map $M_x \rightarrow M_x, \ m \mapsto p\, m$. 

The next observation is pivotal in the present context.

\proclaim{Observation 2.2}  Let $\A$ be a poset, $X$ a topological space, and $f: X \rightarrow  \A$ an upper semicontinuous map whose image is well partially ordered {\rm (meaning that $\Im(f)$ does not contain any infinite strictly descending chain and every nonempty subset has only finitely many minimal elements)}.

Then the pre-images $f^{-1}\bigl(( - \infty, a) \bigr)$ and $f^{-1}\bigl(( - \infty, a] \bigr)$ for $a \in \A$ are open in $X$.  In particular, given any irreducible subset $\U$ of $X$, the restriction of $f$ to $\U$ is generically constant, and the generic value of $f$ on $\,\U$ is 
$$\min\{f(x) \mid x \in \U\}.$$
\endproclaim

\demo{Proof}  We address openness of $f^{-1}\bigl( ( - \infty, a ] \bigr)$.  The case where  $( - \infty, a ]$ contains $\Im(f)$ is trivial.  Otherwise, let $c_1, \dots, c_m$ be the minimal elements of the set 
$\Im(f) \setminus ( - \infty, a]$; the number of such elements is finite by hypothesis.  Then $f^{-1}\bigl( ( - \infty, a ] \bigr)$ is the complement in $\modlad$ of the finite union $\bigcup_{i \le m} f^{-1}\bigl( [c_i, \infty) \bigr)$ of closed sets.

Now let $\U$ be an irreducible subset of $X$ (in particular, $\U \ne \varnothing$), and let $c_1, \dots, c_m$ be the distinct minimal elements in $f(\U)$. Then the sets $\U \cap f^{-1} \bigl( ( - \infty, c_i] \bigr) = f^{-1}(c_i)$ are non-empty and open in $\U$.  Since they are pairwise disjoint, we conclude that $m = 1$.  This makes $c_1$ the generic value of $f$ on $\U$. \qed 
\enddemo

Observation 2.2 entails that, for any noetherian topological space $X$ and any minimal element $a \in \Im(f)$, the closure $\overline{f^{-1}(a)}$ is a finite union of irreducible components of $X$.  This motivates the following terminology.

\definition{Definition 2.3}  Let $X$ be a noetherian topological space, $\A$ a poset, and $f: X \rightarrow \A$ upper semicontinuous.  We say that $f$ {\it detects an irreducible component $\,\C$ of $X$\/} in case the generic value of $f$ on $\C$ is minimal in $\Im(f)$; equivalently, $\C \cap f^{-1}(a) \ne \varnothing$ for some minimal element $a \in \Im(f)$.  Further, $f$ is said to {\it separate irreducible components\/} in case $f^{-1}(a)$ is irreducible for every minimal element $a \in \Im(f)$.  
\enddefinition  

We proceed to the upper semicontinuous maps on $X = \modlad$ which will turn out to detect and separate all irreducible components in case $\la$ is truncated local. 

\proclaim{Observation 2.4} Let $\Seq_{\bd}$ be the set of all semisimple sequences with dimension vector $\bd$, endowed with the dominance order.  Then the map
$$\modlad \longrightarrow \Seq_{\bd}, \ \ \ x \mapsto \SS(M_x)$$
is upper semicontinuous.  Its image coincides with the set of realizable semisimple sequences in $\Seq_{\bd}$.  

Analogously, the map
$$\modlad \longrightarrow \Seq_{\bd}, \ \ \ x \mapsto \SS^*(M_x)$$
is upper semicontinuous, and its image is the set of duals of those semisimple sequences which are realizable in the category of right $\la$-modules.
\endproclaim

\demo{Proof}   It was proved in \cite{\hier, Observation 2.11} that, for any semisimple sequence $\SS$ with $\laySS \ne \varnothing$, the set $\bigcup_{\SS' \ge \SS} \laySSprime$ is closed in $\modlad$.  Upper semi-continuity of socle layerings follows by duality.  The final claim, regarding the image of $\SS^*$, is due to the equality
$$\SS^*(M) = \bigl(D(\SS_0(D(M))), \cdots, D(\SS_L(D(M))) \bigr).  \ \qed$$ 
\enddemo 

We combine the two maps  of Observation 2.4 to an upper semicontinuous map
$$\Theta : \modlad \longrightarrow \Seq_{\bd} \times \Seq_{\bd}, \ \ \ x \mapsto \bigl(\SS(M_x), \SS^*(M_x) \bigr)$$
which we will find to be particularly discerning with regard to irreducible components.  Here, $\Seq_{\bd} \times \Seq_{\bd}$ is equipped with the componentwise dominance order.   In light of Observation 2.2, we wish to identify the pairs $\bigl(\SS, \SS^* \bigr)$ which are minimal in the image of $\Theta$, i.e., minimal in the set
$$ \radsocbd =   \{ \bigl(\SS(M), \SS^*(M) \bigr) \mid M \in \lamod,\, \udim M = \bd \}.$$  
This will permit us to gauge the potential of the map $\Theta$ towards detection and separation of components in the sense of Definition 2.3.  We conclude the section with several lemmas facilitating this task.

\proclaim{Lemma 2.5}  Suppose $M$ is a $\la$-module with dimension vector $\bd$.
\smallskip
\noindent {\rm \bf{(a)}} First suppose that $\SS(M) = (\SS_0, \dots, \SS_L)$.  Then 

$\bullet$  $(\SS_L, \SS_{L-1}, \dots, \SS_0) \le \SS^*(M)$.  

$\bullet$ If $\, \SS^*(M) = (\SS_L, \dots, \SS_0)$, then $\bigl(\SS(M), \SS^*(M) \bigr)$ is minimal in $\radsocbd$.
\smallskip

\noindent {\rm \bf{(b)}} Next suppose that $\SS^*(M) = (\SS^*_0, \dots, \SS^*_L)$.  Then 

$\bullet$  $(\SS^*_L, \SS^*_{L-1}, \dots, \SS^*_0) \le \SS(M)$.  

$\bullet$ If $\, \SS(M) = (\SS^*_L, \dots, \SS^*_0)$, then $\bigl(\SS(M), \SS^*(M)\bigr)$ is minimal in $\radsocbd$.
\endproclaim

\demo{Proof}  We verify (a), part (b) being dual.  A straightforward induction shows that $J^{L - l} M \subseteq \soc_l M$ for $0 \le l \le L$.  (Indeed, this is obvious for $l = 0$, and the inclusion for $l$ implies that the quotient 
$(J^{L- l - 1} M + \soc_l M) / \soc_l M$ 
is annihilated by $J$, whence $J^{L - (l+1)} M \subseteq \soc_{l+1} M$.)  Thus the direct sum of the simple composition factors of $J^{L - l} M$ is contained in the direct sum of the simple composition factors of $\soc_l M$.  Since the former direct sum is 
$$\bigoplus_{0 \le j \le l} J^{L -j} M / J^{L - j + 1} M = \bigoplus_{0 \le j \le l} \SS_{L-j}$$ and the latter equals 
$$\soc_0 M \oplus \bigoplus_{1 \le j \le l} \soc_j M / \soc_{j+1} M = \bigoplus_{0 \le j \le l} \SS_j^*(M),$$ 
the first claim under (a) follows.

To justify the second, suppose $\SS^*(M) =  (\SS_L, \dots, \SS_0)$, and let $N$ be a module with dimension vector $\bd$ such that $\bigl(\SS(N), \SS^*(N)\bigr) \le \bigl( \SS(M), \SS^*(M) \bigr)$.  Due to the inequality of the first entries, $\SS(N) \le \SS(M)$, and the preceding claim for $N$, we obtain $\soc_0 M = J^L M \subseteq J^L N \subseteq \soc_0 N$.  Using the inequality of the second entries, we deduce $\SS^*(M)_0 = \soc_0 M = \soc_0 N = \SS^*(N)_0$, whence also $J^L M = J^L N$.  Analogously, one finds $\soc_1 M / \soc_0 M = J^{L-1} M/ J^L M \subseteq J^{L-1} N / J^L N \subseteq  \soc_1 N/ \soc_0 N$, and a repeat of the initial argument shows that these inclusions are again equalities.  An obvious induction thus yields equality of the pairs $\bigl(\SS(M), \SS^*(M) \bigr)$ and $\bigl(\SS(N), \SS^*(N) \bigr)$.  \qed 
\enddemo 

\proclaim{Corollary 2.6}  Let $\SS = (\SS_0, \dots, \SS_L)$ be any realizable semisimple sequence with dimension vector $\bd$ and $\C$ an irreducible component of $\laySS$.  If $\,\SS^{\text{op}}: = (\SS_L, \dots, \SS_0)$ arises as the socle layering of a module in $\C$, then the closure of $\C$ is an irreducible component of $\modlad$ and $\SS^{\text{op}}$ is the generic socle layering of this component. 
\endproclaim

\demo{Proof}  Let $\D$ be an irreducible component of $\modlad$ containing $\C$.  By Observation 2.2 and Lemma 2.5(a), $\SS^{\text{op}}$ is the generic socle layering of $\D$.  From part (b) we now conclude that $\SS$ is the generic radical layering of $\D$, that is $\D = \overline{\D \cap \laySS} = \overline{\C}$. \qed
\enddemo
 
The upcoming technical points will resurface in the proof of Theorem A.  

\proclaim{Lemma 2.7}  For $M \in \lamod$ the following conditions are equivalent:
\smallskip

{\rm \bf{(a)}} $\SS^*(M) \ne \bigl(\SS_L(M), \dots, \SS_0 (M) \bigr)$.
\smallskip

{\rm \bf{(b)}}  There exist an index $\rho \in \{1, \dots, L\}$ and an element $x \in M \setminus J^\rho M$

\qquad  such that $Jx \subseteq J^{\rho +1} M$.
\smallskip
  
{\rm \bf{(c)}}  There exists an index $\rho \in \{0, \dots, L-1\}$ such that $J^{\rho} M/ J^{\rho + 2} M$ has a simple direct 

\qquad summand.
\endproclaim

\demo{Proof}  Suppose (a) holds, and let $m \ge 0$ be minimal with the property that 
$$\bigl(\SS^*_0(M), \dots, \SS^*_m(M) \bigr) \ne \bigl(\SS_L(M), \dots, \SS_{L-m}(M) \bigr).$$  
Since $J^{L - l} M \subseteq \soc_l M$ for $0 \le l \le L$ by the proof for Lemma 2.5, this amounts to $\soc_j M = J^{L-j} M$ for $j < m$ and $\soc_m M \supsetneqq J^{L-m} M$.  Hence 
$$\soc_m M  / \soc_{m-1} M = \soc \bigl(M / J^{L - (m-1)} M \bigr) \supsetneqq J^{L-m} M / J^{L - m + 1} M,$$
and we infer that (b) is satisfied with $\rho = L - m$.

The implication ``(b)$\implies$(a)" is proved analogously, and the equivalence of (b) and (c) is obvious.  \qed
\enddemo

\proclaim{Lemma 2.8} Suppose that $\SS = (\SS_0, \dots, \SS_L)$ is a realizable semisimple sequence with $\udim \SS = \bd$ such that $\laySS$ is irreducible.  Moreover, suppose that, for each $l \in \{1, \dots, L\}$, the semisimple module $\SS_{l-1}$ embeds into the first socle layer $\SS^*_1(E_l) = \soc_1(E_l)/ \soc_0(E_l)$ of the injective envelope $E_l$ of $\SS_l$.  Then the generic socle layering of the modules in $\laySS$ is $(\SS_L, \dots, \SS_0)$, and $\overline{\laySS}$ is an irreducible component of $\modlad$.  
\endproclaim

\demo{Proof}  In light of Corollary 2.6, the final assertion follows from the claim concerning the generic socle layering of $\laySS$.  By Lemma 2.7, it thus suffices to show that, generically, the modules $M \in \laySS$ satisfy the following condition:  For each $l \in \{1, \dots, L\}$, the subquotient  $J^{l-1} M / J^{l+1} M$ in $\Mod (\SS_{l-1}, \SS_l)$ is free of simple direct summands.  But by hypothesis, $\soc_1 E_l$  contains a submodule with socle layering $(\SS_{l-1}, \SS_l)$, showing that, generically, the modules in $\Mod (\SS_{l-1}, \SS_l)$ have socle $\SS_l$. \qed \enddemo

\head 3.  Narrowing the focus to truncated path algebras \endhead

The purpose of this section is twofold:  Namely, to prepare for applications of Theorem A, and to assemble further background information for its proof.
In the following, we assume $\la = KQ/I$ to be a truncated path algebra with $J^{L+1} = 0$.

\subhead A.  The role of the varieties $\laySS$ and $\biggrassSS$ in the truncated case \endsubhead

In the present situation, the upper semicontinuous map $\Theta$ of Section 2 is known to separate the irreducible components of the module varieties; this is immediate from Theorem 3.1 below.  The initial assertions of this theorem were proved in \cite{\BHT, Theorem 5.3, Corollaries 5.4 and 5.6} by way of the Grassmannian parametrizing varieties; the final statement follows in light of Section 2.
  
\proclaim{Theorem 3.1.  Irreducibility of the $\laySS$} Suppose that $\la$ is a truncated path algebra.
For every realizable semisimple sequence $\SS$ with $\udim \SS  = \bd$, the subvariety $\laySS$ of $\modlad$ is irreducible, unirational and smooth.  Analogous statements hold for the subvariety $\,\biggrassSS$ of $\,\biggrasslad$.  In particular, every irreducible component of $\modlad$ arises as the closure $\,\overline{\laySS}$ for some semisimple sequence $\SS$. 

A sufficient condition for $\, \SS$ to be the generic radical layering of an irreducible component of $\modlad$ is:  \ $\SS = \SS(M)$ such that $\bigl(\SS(M), \SS^*(M) \bigr)$ is minimal in $\radsocbd$.
\qed
\endproclaim

We are thus prompted to explore equivalent conditions for $\laySS$ to close off to an irreducible component of $\modlad$.

The first entries of the pairs in $\radsocbd$, namely the realizable semisimple sequences, can be listed by mere inspection of the quiver in the truncated scenario.  (The use of skeleta  --  see Subsection B  --  meets the need to keep track of counts.)  Indeed, the following observation regarding realizability is straightforward. 

\proclaim{Observation 3.2}  Let $\la$ be a truncated path algebra, $\SS = (\SS_0, \dots, \SS_L)$ a semisimple sequence in $\lamod$, and $P_l$ a projective cover of $\SS_l$ for $0 \le l \le L$. Then the following conditions are equivalent:

$\bullet$ $\SS$ is realizable, that is, $\SS = \SS(M)$ for some $M \in \lamod$.

$\bullet$ For each $l \in \{1, \dots, L\}$, the sequence $(\SS_{l-1}, \SS_l)$ is realizable over the algebra $\la / J^2$.

$\bullet$  For $l \in \{1, \dots, L\}$, the entry $\SS_l$ is contained in the first radical layer $J P_{l-1} / J^2 P_{l-1}$ of $P_{l-1}$, that is:  If $\SS_l = \bigoplus_{1 \le i \le n} S_i^{m_{li}}$, then 
$$m_{lj} \ \ \le \ \ \sum_{1 \le i \le n} m_{l-1, i} \cdot \#(\text{arrows\ } e_i \rightarrow e_j) \ \ \ \ \text{for} \ 1 \le j \le n. \ \  \qed$$
\endproclaim 

The two final conditions listed in Observation 3.2 are always (beyond the truncated case) equivalent, and they are always necessary for realizabilty of $\SS$.  However, they need not be sufficient for more general types of algebras.  We contrast Theorem 3.1 with the situation of a non-truncated algebra.

\definition{Example 3.3} Consider the algebra $\la = KQ/\langle \beta_i \alpha_j \mid i \ne j \rangle$, where $Q$ is the quiver with three vertices $e_1, e_2, e_3$, two arrows, $\alpha_1, \alpha_2$, from $e_1$ to $e_2$, and two arrows, $\beta_1, \beta_2$, from $e_2$ to $e_3$.  If $\SS = (S_1, S_2, S_3)$, then $\laySS$ is reducible.  Indeed, this variety has two irreducible components, namely the orbit closures of the uniserial modules with the following graphs:
$$\xymatrixrowsep{1.5pc}\xymatrixcolsep{2pc}
\xymatrix{
1 \edge[d]^{\alpha_1} && & &&1 \edge[d]^{\alpha_2}  \\
2 \edge[d]^{\beta_1} && &\txt{and} &&2 \edge[d]^{\beta_2}  \\
3 && & &&3
}$$
\noindent Moreover, the sequence $\SS' = (S_1, S_2, S_3^2)$ fails to be realizable, even though it satisfies the final conditions of Observation 3.2. \qed
\enddefinition
\smallskip

\subhead B.  Skeleta and generic projective presentations over truncated path algebras  \endsubhead

We recall some concepts which were defined for and applied to arbitrary basic finite dimensional algebras (see, e.g., \cite{\BHT} or \cite{\hier}).  However, since restriction to the truncated case simplifies the picture, we pare down the general definitions and results to take advantage of the current situation.    

Very roughly, a skeleton of a $\la$-module $M$ is a path basis of $M$ which lives in a fixed projective cover $P$ of $M$, displays the radical layering of $M$, and is ``shared" by all modules in a dense open subset of the irreducible variety $\Mod\, \SS(M)$.  So, the set of skeleta of a module is a generic attribute of the modules in any irreducible component of $\modlad$.

\definition{Definition 3.4 and comments}  Let $\la = KQ/ \langle \text{the paths of length \ } L+1 \rangle$ be a truncated path algebra and $\SS = (\SS_0, \dots, \SS_L)$ a semisimple sequence with dimension vector $\bd$ and $\dim \SS_0 = t$.  Fix a projective cover $P = \bigoplus_{1 \le r \le t} \la z_r$ of $\,\SS_0$, where $z_1, \dots, z_t$ is a full sequence of top elements of $P$.  (We point out that, typically,  $P$ is a proper direct summand of the projective module $\bP$ on which we based the definition of $\biggrass_\bd(\la)$ in Section 1; note that $\udim \bP / J\bP = \bd$.) Any {\it nonzero} element of the form $p z_r \in P$, where $p$ is a path in $Q$, is called a {\it path in\/} $P$, and we set $\len(p z_r) = \len(p)$.  Note that this definition is unambiguous, since $\la$ is graded by path lengths.  In particular, the paths of length $0$ in $P$ are precisely those of the form $z_r = e(r) z_r$, $1 \le r \le t$, where $e(r)$ is the primitive idempotent that norms $z_r$.  The {\it endpoint\/}, end($p z_r$), of a path $p z_r$ is that of $p$. 

{\bf (a)} An ({\it abstract\/}) {\it skeleton with radical layering\/} $\SS$ is a set $\S$ of paths in $P$ such that

\quad$\bullet$  $\S$ is closed under initial subpaths, meaning:  whenever $p = p_2 p_1$ for paths $p_i$ in $Q$ such that $p_2 p_1 z_r \in \S$, it follows that $p_1 z_r \in \S$.

\quad$\bullet$  For any $l \in \{0, \dots, L\}$ and $i \in \{1, \dots, n\}$, the number of those paths of length $l$ in $\S$ which end in $e_i$ equals the multiplicity of $S_i$ in $\SS_l$.

{\bf (b)}  A {\it skeleton\/} of $M \in \lamod$ is a skeleton $\S$ with radical layering $\SS(M)$ such that $M \cong P/C$ and $\{pz_r + C \mid p z_r \in \S\}$ is a basis for $P/C$.
\enddefinition

Note that the set of abstract skeleta with any fixed radical layering is finite and that every $\la$-module has at least one skeleton.  Conversely, by Observation 3.2, each abstract skeleton arises as the skeleton of a module.  In light of the graphical rendering of skeleta described below, they provide visually suggestive labels for the realizable semisimple sequences over $\la$: Indeed, $\SS$ is realizable if and only if there is an abstract skeleton with radical layering $\SS$.  

However, the primary motivation for considering skeleta lies in the fact that suitable subvarieties $\biggrassS$ of $\biggrassSS$, consisting of points that correspond to modules with skeleton $\S$, constitute a particularly useful open affine cover of $\biggrassSS$;  the charts $\biggrassS$ of this cover are always stable under the action of the unipotent radical of $\aut_\la(\bP)$ and serve to translate geometric information regarding $\biggrassSS$ into algebraic information on the corresponding representations. In the case of a truncated path algebra, they are affine spaces which exhibit generic minimal projective presentations of the modules in $\biggrassSS$ (resp., $\laySS$).  The connection between skeleta and minimal projective presentations will be spelled out in Observation 3.6 below and used in Supplement B to Theorem A for instance (Section 4).  For further detail, in wider generality, we refer to \cite{\BHT} and \cite{\hier}. 

It is convenient to identify skeleta $\S$ with labeled and layered graphs.  The latter not only convey the radical layering of any module $M$ with skeleton $\S$ at a glance, but also the preferred basis:  The graph of $\S$ is a forest consisting of trees $\T_1, \dots, \T_t$, the tree $\T_r$ having root $z_r$.  The $\T_r$ are (informally) defined by the following property: For each $r \le t$, the edge paths of length $\ge 0$ starting in the root $z_r$ of $\T_r$  are precisely the paths of the form $p z_r$ in $\S$.  (The first of the defining properties of a skeleton guarantees that, for any $r \le t$, the paths $p z_r \in \S$ indeed assemble to a tree.)

We illustrate these concepts in

\definition{Example 3.5}  Let $\la$ be the truncated path algebra of Loewy length $4$ over $\CC$ based on the quiver
$$\xymatrixrowsep{1pc}\xymatrixcolsep{2pc}
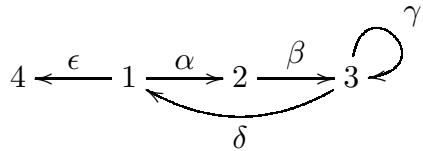
\xymatrix{
4 &1 \ar[l]_{\epsilon} \ar[r]^{\alpha} &2 \ar[r]^{\beta} &3 \ar@(u,r)^{\gamma} \ar@/^1pc/[ll]^{\delta}
}$$
\noindent and let $\SS = (S_1 \oplus S_2, S_2 \oplus S_3 \oplus S_4, S_3^2, S_1 \oplus S_3)$.  There are precisely $4$ skeleta $\S$ with radical layering $\SS$, namely
$$\xymatrixrowsep{1.5pc}\xymatrixcolsep{0.25pc}
\xymatrix{
1 \edge[d]_{\alpha} \edge[dr]^{\epsilon} &&2 \edge[d]^{\beta}  &&&&&&1 \edge[d]_{\alpha} \edge[dr]^{\epsilon} &&2 \edge[d]^{\beta}  &&&&&&1 \edge[d]_{\alpha} \edge[dr]^{\epsilon} &&2 \edge[d]^{\beta}  &&&&&&1 \edge[d]_{\alpha} \edge[dr]^{\epsilon} &&2 \edge[d]^{\beta}  \\
2 \edge[d]_{\beta} &4 &3 \edge[d]^{\gamma}  &&&&&&2 \edge[d]_{\beta} &4 &3 \edge[d]^{\gamma}  &&&&&&2 \edge[d]_{\beta} &4 &3 \edge[d]^{\gamma}  &&&&&&2 \edge[d]_{\beta} &4 &3 \edge[d]^{\gamma}  \\
3 \edge[d]_{\delta} \edge[dr]^{\gamma} &&3  &&&&&&3 \edge[d]_{\delta} &&3 \edge[d]^{\gamma}  &&&&&&3 \edge[d]_{\gamma} &&3 \edge[d]^{\delta}  &&&&&&3 &&3 \edge[d]_{\delta} \edge[dr]^{\gamma}  \\
1 &3 &  &&&&&&1 &&3  &&&&&&3 &&1  &&& &&&&&1 &3
}$$
\noindent  The generic socle layering of the modules in $\laySS$ is $\SS^* = (S_1 \oplus S_3^2 \oplus S_4, S_3^2, S_2^2, S_1)$, and the pair $(\SS, \SS^*)$ is a minimal element of $\radsocbd$, where $\bd = \udim \SS$.  In particular, $\C= \overline{\laySS}$ is an irreducible component of $\modlad$, and generically the modules in $\C$ have all four of the displayed skeleta simultaneously.  The latter follows from the fact that the $\biggrassS$ are open and dense in the irreducible variety $\biggrassSS$ by the preceding remarks and Proposition 1.2. \qed
\enddefinition 

The upcoming observation is covered by \cite{\BHT, Section 3 and Theorem 5.3}.  We continue to assume that $\la$ is a truncated path algebra.

\proclaim{Observation 3.6.  Minimal projective presentations of the modules with skeleton $\S$}  Let $M$ be a $\la$-module with radical layering $\SS$, and suppose that $\S \subseteq P$ is a skeleton of $M$; here $P = \bigoplus_{1 \le r \le t} \la z_r$ is the fixed projective cover of $\SS_0$ as in Definition {\rm 3.4}, which makes $P$ a projective cover also of $M$.  We say that a path $u z_r \in P$ is $\S$-\underbar{critical} if $u z_r \notin \S$, while all proper initial subpaths of $u z_r$ belong to $\S$.  For any $\S$-critical path $u z_r$, let 
$$\S(u z_r) = \{q z_s \in \S \mid \len(q z_s) \ge \len(u z_r) \text{\ and\ } \operatorname{end}(u z_r) = \operatorname{end}(q z_s)\}.$$ 
Then there is a unique family of scalars, $c = \bigl(c_{u z_r, q z_s}\bigr)$, where $uz_r$ traces the $\S$-critical paths and $q z_s$ traces $\S(u z_r)$, respectively, such that 
$M \cong P/U(c)$, where $U(c)$ is the submodule of $P$ which is generated by the differences
$$u z_r\ \  - \  \ \sum_{q z_s \in \S(u z_r)} c_{u z_r, q z_s}\, q z_s \ \ \ \ \text{for} \ \ uz_r \ \, \S\text{-critical}.$$
Conversely, for any choice $c$ of scalars, $P/U(c)$ is a module with skeleton $\S$.  \qed
\endproclaim

For a proof of the upcoming Theorem 3.7  --  it motivates our use of the term {\it generic module\/}  --   we refer to \cite{\BHT, Theorems 4.3 and 5.12}.  

\proclaim{Definition of generic modules and Theorem 3.7}  {\rm Let $\la$ be a truncated path algebra, $\SS = (\SS_0, \dots, \SS_L)$ a realizable semisimple sequence, and $P = \bigoplus_{1 \le r \le t} \la z_r$ the projective cover of $\SS_0$ as in Definition {\rm 3.4}.  Moreover, given any skeleton $\S$ with radical layering $\SS$, let $\Ktilde$ be an algebraic closure of $K(X)$, where $X = \bigl(X_{u z_r, q z_s})$ is a family of independent variables over $K$ with indices ranging over the set of pairs 
$$N : = \{(u z_r, q z_s) \mid u z_r \in P\ \text{is}\  \S\text{-critical and}\  q z_s \in \S(u z_r)\}.$$  
We will refer to the following module $G(X)$ over $\latilde := \Ktilde \otimes_K \la$ as a {\it generic module\/} for $\laySS$: namely, $G(X) := \bigl( \Ktilde \otimes_K P \bigr) / U(X)$, where $U(X)$ is the $\latilde$-submodule of $\Ktilde \otimes_K P$ which is generated by the differences
$$u z_r\ \  - \  \ \sum_{q z_s \in \S(u z_r)} X_{u z_r, q z_s}\ q z_s.$$

\noindent To motivate this terminology, embed $\lamod$ into $\latilde\text{-}\operatorname{mod}$ via $M \mapsto \Ktilde \otimes_K M$.}  Then the module $G(X)$ has all those generic properties of the modules in $\laySS$ which are preserved by Morita self-equivalences of $\latilde\text{-}\operatorname{mod}$; more strongly, this extends to all generic properties which are preserved by Morita self-equivalences induced by field automorphisms in $\operatorname{Gal}\, (\Ktilde/K)$.  Moreover, with respect to this condition, $G(X)$ is unique up to a Morita self-equivalence of the restricted type.  \qed
\endproclaim

Note that the fields $K$ and $\Ktilde$ are isomorphic in case $K$ has infinite transcendence degree over its prime field.

\definition{Return to Example 3.5}  A generic module for $\laySS$ is $G(X) = (\latilde e_1 \oplus \latilde e_2)/U(X)$, where $U(X)$ is generated by $\delta \beta - X_1 \delta \beta \alpha$, $\delta \gamma \beta - X_2 \delta \beta \alpha$, and $\gamma^2 \beta - X_3 \gamma \beta \alpha$.  The presentation of the module $G(X)$ can be simplified, but $G(X)$ is indecomposable.  Therefore the modules in the irreducible component $\C$ of $\modlad$ are generically indecomposable. \qed
\enddefinition

We still refer to $G(X)$ as {\it generic for\/} $\laySS$ if the algebraically closed base field $\Ktilde$ of $\Ktilde\otimes_K \la$ has higher transcendence degree over $K$ than is called for by the preceding definition.  This comment bears on the next statement.

\proclaim{Proposition 3.8}  Retain the hypotheses and notation of Theorem {\rm 3.7}, and denote by $J$ the Jacobson radical also of the $\Ktilde$-algebra $\Ktilde \otimes_K \la$.
If $G(X)$ is a generic module for $\laySS$, then each of the subfactors $J^{\rho} G(X) / J^{\tau + 1} G(X)$,  for $0 \le \rho < \tau \le L$, is a generic module for the irreducible  variety $\Mod (\SS_{\rho}, \dots,  \SS_{\tau})$ over $\la/J^{\tau - \rho + 1}$.  

In particular, generic decomposability of the $\la$-modules in $\laySS$ entails generic decomposability of the modules in $\Mod (\SS_{\rho}, \dots,  \SS_{\tau + 1})$. \qed
\endproclaim

\head 4. The local case.  Consequences of Theorem A and examples \endhead

Barring the final example of this section, $\la$ denotes the local truncated path algebra of Loewy length $L+1$ over $K$ with $\dim_K J/J^2 = r$; i.e.,  $\la$ is based on the following quiver:

$$\xymatrixrowsep{0.75pc}\xymatrixcolsep{0.3pc}
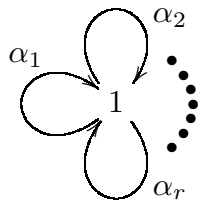
\xymatrix{
 & \ar@{{}{*}}@/^1pc/[dd] \\
1 \ar@'{@+{[0,0]+(-6,-6)} @+{[0,0]+(-15,0)} @+{[0,0]+(-6,6)}}^(0.7){\alpha_1}
\ar@'{@+{[0,0]+(-6,6)} @+{[0,0]+(0,15)} @+{[0,0]+(6,6)}}^(0.7){\alpha_2}
\ar@'{@+{[0,0]+(6,-6)} @+{[0,0]+(0,-15)} @+{[0,0]+(-6,-6)}}^(0.3){\alpha_r}  \\
 &
}$$

\noindent  Throughout, we assume the number $r$ of loops to be at least $2$.
Theorem A was stated in the introduction.

In light of Section 3, this theorem provides us with generic minimal projective presentations of the modules in the irreducible components of the varieties $\modlasd$.  In stating our claim, we refer to the notation of Observation 3.6.

\proclaim{Supplement B to Theorem A}  Let $\SS = (\SS_0, \dots, \SS_L)$ be a $d$-dimensional semisimple sequence with $\dim \SS_l \le r \cdot \dim \SS_{l-1}$ and $\,\dim \SS_{l-1} \le r \cdot \dim \SS_l$ for $1 \le l \le L$, and let $\S$ be any skeleton with radical layering $\SS$.  

Generically, the modules in the irreducible component $\C = \overline{\laySS}$ of $\modlasd$ have a minimal projective presentation of the form $P/U(c)$, where $P$ is a projective cover of $\, \SS_0$ and $U(c)$ is as described in Observation {\rm 3.6}. In fact, there is a dense open subset $\,\U$ of $\,\C$ such that the modules in $\,\U$ are precisely those isomorphic to $P/U(c)$, where $c = (c_{u z_r, q z_s})$ traces the affine $K$-space of dimension $N = | \{ (uz_r , qz_s) \mid uz_r \text{\rm{\  a \, }} \S\text{\rm{-critical path,\ }} q z_s \in \S(u z_r)\}|$. 
\endproclaim

\demo{Proof}   That $\C$ is an irreducible component of $\modlasd$ follows from Theorem A.  Moreover, $\laySS$ is locally closed, and hence is open in its closure.  Consequently, Proposition 1.2 guarantees that the subvariety $\biggrassSS$ of the projective variety $\biggrass_d(\la)$ closes off to an irreducible component $\D$ of the latter and is open in its closure.  Combining \cite{\hier,  Theorem 3.17} with \cite{\BHT, Theorem 5.3}, we find that $\biggrassSS$ contains a dense open subset $\biggrassS$ which consists of points corresponding to modules with skeleton $\S$. The subvariety $\biggrassS$ in turn contains a subvariety $\grassS$ with the following properties: $\bullet$  $\grassS$ is an affine variety isomorphic to $\AA^N$, with a coordinate system indexed by the set $N =  \{ (uz_r , qz_s) \mid uz_r \text{\rm{\ a \ }} \S\text{\rm{-critical path,\ }} q z_s \in \S(u z_r)\}$, and $\bullet$ the module corresponding to any point $c = (c_{u z_r, q z_s}) \in \grassS$ is $P/U(c)$ (see \cite{\hier, Section 3}).  Let $\V$ be the closure of $\biggrassS$ under the $\aut_\la(\bP)$-action.  Then $\V$ is in turn open in $\biggrassSS$, and thus in $\D$. Another application of Proposition 1.2 now shows the subvariety $\U = \Phi(\V)$ of $\C$ to be as claimed.  (Note that $\biggrassS$ is not closed under the $\aut_\la(\bP)$-action of $\biggrass_d(\la)$ in general, whence Proposition 1.2 does not permit a straightforward shift of information from the projective to the affine setting.)\qed  
\enddemo

The next corollary to Theorem A is evident (separately deal with the cases $d \le L+1$ and $d > L+1$).  The subsequences consisting of two consecutive entries $\SS_{l-1}, \SS_l$ will play a role in deciding generic indecomposability of the modules in $\laySS$.    

\proclaim{Corollary C}  Let $\la$ be a local truncated path algebra with Loewy length $L+1$.  A semisimple sequence $\,\SS = (\SS_0, \dots, \SS_L)$ is the generic radical layering of an irreducible component of $\modlasd$ if and only if, for each $\l \in \{1, \dots, L\}$, the two-term sequence $(\SS_{l-1}, \SS_l)$ is the generic radical layering of an irreducible component $\C_l$ of $\, \Mod_{d^{(l)}} \la / J^2$, where $d^{(l)} = \dim \SS_{l-1} + \dim \SS_l$.  
\qed
\endproclaim

Generic modules for the components $\C_l$ are exhibited in Proposition 3.8.  

Clearly, Corollary C extends to arbitrary segments of $\SS$:  Namely, $\SS$ is the generic radical layering of an irreducible component of $\modlasd$ if and only if, for any choice of $\rho < \tau$ in $\{0, \dots, L\}$, the trimmed sequence $(\SS_\rho, \dots, \SS_{\tau})$ is the generic radical layering of an irreducible component of the variety of $(\la/J^{\tau - \rho + 1})$-modules of dimension $\sum_{\rho \le l \le \tau} \dim \SS_l$.  

\proclaim{Corollary D}  Again, let $\la$ be a local truncated path algebra with Loewy length $L+1$.  Moreover, suppose $d > L+1$, and let $\C$ be an irreducible component of $\modlasd$.  Then the modules $M$ in $\C$ generically satisfy the following condition:  $J^L x \ne 0$ whenever $x \in M \setminus JM$.

In particular: There exists a dense open subset $\,\U \subseteq \C$ such that every module in $\U$ has a decomposition  into indecomposable direct summands \underbar{all} of which have Loewy length $L+1$. \endproclaim

\demo{Proof}  Suppose $\C = \overline{\laySS}$, and let $M$ be a module in $\laySS$ whose socle layering is the generic one, i.e., $\SS^*(M) = \bigl(\SS_L(M), \dots, \SS_0(M) \bigr)$ by Theorem A.  Moreover, suppose $x \in M \setminus JM$.  Let $s$ be minimal with $J^s x = 0$, and $t$ maximal with $J^{s-1}x \subseteq J^t M$.  Then $L+1 \ge s > t \ge 0$, and $J^t M/ J^{t+2} M$ has a nontrivial semisimple direct summand, contained in the canonical image of $J^{s-1} x$.  Hence we conclude from Lemma 2.7 that $t=L$ and $s = L+1$.  \qed 
\enddemo

It is easy to see that Corollary D has no analogue for general truncated path algebras. Corollaries C and D, in turn, have an interesting consequence. 
Namely, in many cases, the test for generic indecomposability of the modules in an irreducible component of $\modlasd$ may be played back to the generalized Kronecker quiver $\Qhat$:
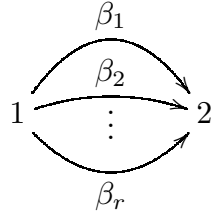
$$\xymatrixrowsep{1pc}\xymatrixcolsep{4pc}
\xymatrix{
1 \ar@/^2pc/[r]^{\beta_1} \ar@/^/[r]^{\beta_2} \ar@{}[r]|{\vdots} \ar@/_1.5pc/[r]_{\beta_r}  &2
}$$
\noindent The Schur roots of $\Qhat$, i.e., the dimension vectors $(d_1,d_2)$ with the property that the modules in $\Mod_{(d_1, d_2)} K \Qhat$ are generically indecomposable, were already determined by Kac in \cite{\KacI, Section 2.6}.  
Observe that, whenever $d \le L+1$, the modules in the irreducible variety $\modlasd$ are generically indecomposable by Theorem A.  Beyond that we obtain:

\proclaim{Corollary E} Let $\la$ be as in the preceding corollaries.  Suppose $d > L+1$, and let $\SS$ be the generic radical layering of the modules in an irreducible component $\, \C$ of $\modlasd$.  A sufficient condition for the modules in $\C$ to be generically indecomposable is as follows:  There exists an index $l \in \{1, \dots, L\}$ such that the pair $(\dim \SS_{l-1}, \dim \SS_l)$ is a Schur root of the generalized Kronecker quiver $\Qhat$.

Moreover:  In case at least one of the varieties $\Mod_{(\dim \SS_{l-1}, \dim \SS_l)} K \Qhat$ contains infinitely many orbits of maximal dimension under the action of the pertinent general linear group, then $\, \C$ contains infinitely many $\GL_d$-orbits of maximal dimension.
\endproclaim

\demo{Proof} Suppose that the modules in $\modlasd$ are generically decomposable. In other words, there exists a dense open subset $\U \subseteq \C$ such that every module $M$ in $\,\U$ decomposes in the form $M = M_1 \oplus M_2$ with $\dim M _i$ fixed and nonzero.   Generically, the modules $M$ in $\U$ then satisfy $J^L M_i \ne 0$ for $i = 1,2$ by Corollary D.  In particular, we find that, for each $l \in \{1, \dots, L\}$, the quotient $J^{l-1} M / J^{l+1} M$ is in turn decomposable.  By Corollary C and Proposition 3.8, these quotients trace, for each $l$, a dense open subset of an irreducible component $\C_l$ of $\Mod_{d^{(l)}}(\la / J^2)$, where $d^{(l)} = \dim \SS_{l-1} + \dim \SS_l$.  Now \cite{\BCH, Theorem 5.6(b)} guarantees that the modules in $\Mod_{(\dim \SS_{l-1}, \dim \SS_l)} K\Qhat$ are decomposable for $1 \le l \le L$; in other words, none of the dimension vectors $(\dim \SS_{l-1}, \dim \SS_l)$ is a Schur root of $\Qhat$. 

The final assertion follows from \cite{\BCH, Theorem 5.6(c)} by analogous reasoning.  \qed 
\enddemo  

The sufficient condition for generic indecomposability given in Corollary E fails to be necessary in general.  Failure may even occur for $r = 2$:

\definition{Example 4.1}  Let $\la$ be the local truncated $\CC$-algebra of Loewy length $3$ based on the quiver with $2$ loops, and denote by $S$ the unique simple in $\lamod$.  For $d = 6$, the semisimple sequence $\SS = (S^2, S^2, S^2)$ is the generic radical layering of an irreducible component of $\modlasd$, but $(2,2)$ fails to be a Schur root of the Kronecker algebra with two arrows.  On the other hand, the modules in $\overline{\laySS}$ are generically indecomposable; indeed, one checks that the endomorphism ring of the generic module $G(X)$ for $\laySS$ has top $\CC$ (see Theorem 3.7).  Given that only finitely many variables can be eliminated from the presentation of $G(X)$, the component $\overline{\laySS}$ moreover contains infinitely many orbits of maximal dimension, a fact that can alternately be gleaned from the Kronecker algebra by way of the final statement of Corollary E.  \qed
\enddefinition

We illustrate Theorem A and Corollary E.  

\definition{Example 4.2}  Let $\la$ be the local truncated path algebra with $r = 3 = L+1$, and $d = 10$.  Then $\modlasd$ has precisely 17 irreducible components.  Indeed, the eligible generic radical layerings are readily listed via Criterion (2) of Theorem A. They are displayed below, the bullets indicating the dimensions of the layers.
$$\xymatrixrowsep{0.1pc}\xymatrixcolsep{0.001pc}
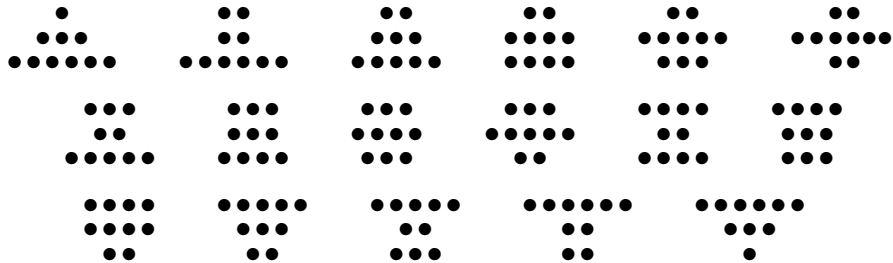
\xymatrix{
 &&\mbull &&& &&&& &&\dbull &\dbull && &&&& &\mbull &\mbull && &&&& &\dbull &\dbull & &&&& &\mbull &\mbull && &&&& &&\dbull &\dbull  \\
 &\mbull &\mbull &\mbull && &&&& &&\dbull &\dbull && &&&& &\dbull &\dbull &\dbull & &&&&\dbull &\dbull &\dbull &\dbull &&&&\dbull &\dbull &\dbull &\dbull &\dbull &&&&\dbull &\dbull &\dbull &\dbull &\dbull &\dbull  \\
\dbull &\dbull &\dbull &\dbull &\dbull &\dbull &&&&\dbull &\dbull &\dbull &\dbull &\dbull &\dbull &&&&\dbull &\dbull &\dbull &\dbull &\dbull &&&&\dbull &\dbull &\dbull &\dbull &&&& &\dbull &\dbull &\dbull & &&&& &&\dbull &\dbull  \\  \\
&&& &\dbull &\dbull &\dbull & &&&&\mbull &\mbull &\mbull & &&&&\mbull &\mbull &\mbull & &&&& &\dbull &\dbull &\dbull & &&&&\dbull &\dbull &\dbull &\dbull &&&&\dbull &\dbull &\dbull &\dbull  \\
&&& &\mbull &\mbull && &&&&\mbull &\mbull &\mbull & &&&&\dbull &\dbull &\dbull &\dbull &&&&\dbull &\dbull &\dbull &\dbull &\dbull &&&& &\dbull &\dbull & &&&&\mbull &\mbull &\mbull &  \\
&&& \dbull &\dbull &\dbull &\dbull &\dbull &&&&\dbull &\dbull &\dbull &\dbull &&&&\mbull &\mbull &\mbull & &&&& &\mbull &\mbull && &&&&\dbull &\dbull &\dbull &\dbull &&&&\mbull &\mbull &\mbull &  \\  \\
&&&& \dbull &\dbull &\dbull &\dbull &&&&\dbull &\dbull &\dbull &\dbull &\dbull &&&&\dbull &\dbull &\dbull &\dbull &\dbull &&&&\dbull &\dbull &\dbull &\dbull &\dbull &\dbull &&&&\dbull &\dbull &\dbull &\dbull &\dbull &\dbull  \\
&&&& \dbull &\dbull &\dbull &\dbull &&&& &\dbull &\dbull &\dbull & &&&& &\mbull &\mbull && &&&& &&\dbull &\dbull && &&&& &\mbull &\mbull &\mbull &  \\
&&&& &\dbull &\dbull & &&&& &\mbull &\mbull && &&&& &\dbull &\dbull &\dbull & &&&& &&\dbull &\dbull && &&&& &&\mbull &
}$$
 All of the components parametrize generically indecomposable modules and have infinitely many $\GL_d$-orbits of maximal dimension; and in all but one of the cases (the last radical layering in the first row), reference to the Kronecker algebra with $3$ arrows yields these findings.   For instance, the second of the irreducible components of $\Mod_{10}(\la)$ is the closure $\, \C$ of $\laySS$, where $\SS = (S^2, S^2, S^6)$.  In this case $(\dim \SS_0, \dim \SS_1) = (2,2)$ is a Schur root of the Kronecker quiver with three arrows;  see \cite{\KacI, Theorem 4(a)}.  That $\C$ does not contain a dense orbit can be gleaned from \cite{\BCH, Theorem 5.6 (c)}. 
 \qed \enddefinition

Next we provide simple instances of the two different deformation strategies underlying the proof of Theorem A. They facilitate visualization of the argument across technical hurdles.

\definition{Example 4.3 illustrating the proof of Lemma 5.2}  Let $\la$ be the local algebra with $r = 3 = L+1$.  Any $\la$-module $M$ with a graph as shown on the left below deforms to a module $\Mtilde$ with a graph as shown on the right.  By this we mean: $M$ belongs to the closure $\overline{\Mod \SS(\Mtilde)}$.
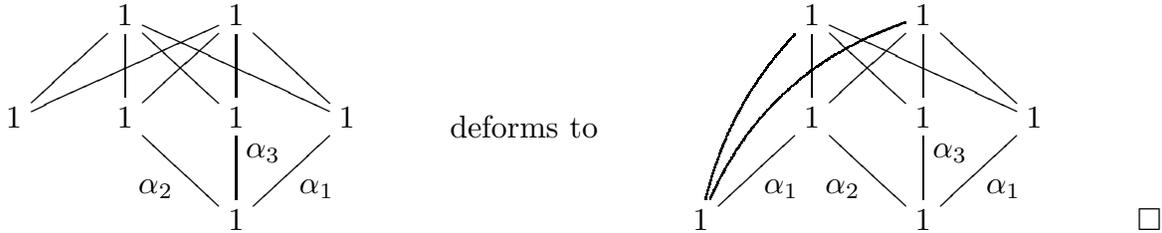
$$\xymatrixrowsep{1.5pc}\xymatrixcolsep{2pc}
\xymatrix{
 &1 \edge[dl] \edge[d] \edge[dr] \edge[drr] &1 \edge[dll] \edge[dl] \edge[d] \edge[dr]  & && &1 \edge@/_/[ddl] \edge[d] \edge[dr] \edge[drr] &1 \edge@/_1pc/[ddll] \edge[dl] \edge[d] \edge[dr]  \\
1 &1 \edge[dr]_{\alpha_2} &1 \edge[d]^(0.33){\alpha_3} &1 \edge[dl]^{\alpha_1}  &\txt{deforms to} &&1 \edge[dl]^{\alpha_1} \edge[dr]_{\alpha_2} &1 \edge[d]^(0.33){\alpha_3} &1 \edge[dl]^{\alpha_1}  \\
 &&1 & &&1 &&1 &&\square
 }$$
\enddefinition

\definition{Example 4.4 illustrating the first step of the proof of Theorem A}  Let $\la$ be the local algebra with $r = 3 = L$.  Any $\la$-module $M$ with a graph as shown on the left below deforms to a $\la$-module $\Mtilde$ with a graph as shown on the right.
$$\xymatrixrowsep{1.5pc}\xymatrixcolsep{2pc}
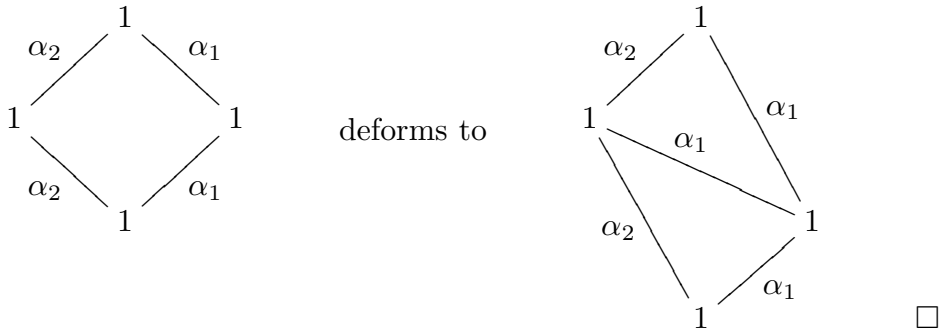
\xymatrix{
 &1 \edge[dl] _{\alpha_2} \edge[dr]^{\alpha_1} & && &1 \edge[dl]_{\alpha_2} \edge[ddr]^{\alpha_1}  \\
1 \edge[dr]_{\alpha_2} &&1 \edge[dl]^{\alpha_1}  &\txt{deforms to} &1 \edge[ddr]_{\alpha_2} \edge[drr]^(0.4){\alpha_1}  \\
 &1 & && &&1 \edge[dl]^{\alpha_1}  \\
 && && &1 &&\square
 }$$
\enddefinition

As mentioned in the introduction, for nonlocal truncated path algebras of Loewy length $\ge 4$, the upper semicontinuous map $(\SS, \SS^*)$ may fall short of detecting all irreducible components of the module varieties.

\definition{Example 4.5}  Let $\la = KQ / \langle \text{all paths of length\ } 4 \rangle$, where $Q$ is the quiver 
$$\xymatrixrowsep{1pc}\xymatrixcolsep{2pc}
\xymatrix{
1 \ar[r]^{\alpha} &2 \ar@/^/[r]^{\beta} &3 \ar@/^/[l]^{\delta} \ar[r]^{\gamma} &4
}$$
For $\bd = (1,1,1,1)$, the variety $\modlad$ consists of two irreducible components cut out by the following generic radical layerings: $\SS  =  (S_1, S_2, S_3, S_4)$ and $\SStilde  = (S_1 \oplus S_3, S_2 \oplus S_4, 0 , 0)$.  Both of the varieties $\laySS$ and $\Mod \SStilde$ contain dense orbits, namely those corresponding to the modules $G$ and $\Gtilde$ determined by the following graphs, respectively; indeed, both $G$ and $\Gtilde$ have only trivial  self-extensions, whence their $\GL_{\bd}$-orbits are open in $\modlad$. 
$$\xymatrixrowsep{1.5pc}\xymatrixcolsep{2pc}
\xymatrix{
1 \edge[d]_{\alpha}   &&1 \edge[d]_{\alpha} &3 \edge[d]^{\gamma}  \\
2 \edge[d]_{\beta}  &&2 \edge[ur]^(0.6){\delta} &4  \\
3 \edge[d]_{\gamma}  \\
4
}$$
\noindent   Therefore the generic socle layerings of $\laySS$ and $\Mod\SStilde$ are
$$\SS^* =  \SS^*(G) = (S_4, S_3, S_2, S_1) \ \ \text{and}\ \ \SStilde^* = \SS^*(\Gtilde) =   (S_2 \oplus S_4, S_1 \oplus S_3, 0, 0),$$
respectively, which shows $(\SS, \SS^*) < (\SStilde, \SStilde^*)$.  On the other hand, $\delta \cdot G = 0$ while $\delta \cdot \Gtilde \ne 0$, whence the corresponding generic triples 
$$\bigl(\SS(G), \SS^*(G), \nullity_\delta(G) \bigr) \ \ \  \text{and} \ \ \ \bigl(\SS(\Gtilde), \SS^*(\Gtilde), \nullity_\delta(\Gtilde) \bigr)$$ 
are not comparable; here $\nullity_\delta G = \dim \ann_G (\delta)$ as in the remarks preceding Observation 2.2.  Thus, the radical-socle layering by itself fails to detect the component $\overline{\Mod\SStilde}$. \qed
\enddefinition

\head 5.  Proof of Theorem A \endhead

Throughout this section, let $\la$ be as in the hypothesis of Theorem A, that is, $\la = KQ / \langle \text{all paths of length}\ L+1 \rangle$, where $Q$ has a single vertex and $Q_1$ consists of $r \ge 2$ loops, say $Q_1 = \{\alpha_1, \dots, \alpha_r\}$.  
The unique (up to isomorphism) simple left $\la$-module is denoted by $S$. 
Further,  $Q_{>0}$ will stand for the set of paths of positive length in $Q$.  
Given any family $(f_\alpha)_{\alpha \in Q_1}$ of $K$-endomorphisms of $K^d$, the following notation will be convenient: whenever $p = \alpha_{i_l} \cdots \alpha_{i_1}$ is a path in $Q_{>0}$, we let $f_p$ be the corresponding composition of maps $f_{\alpha_j}$.

We will repeatedly use the specialization of Observation 3.2 to the local case:  Namely, a semisimple sequence $\SS = (\SS_0, \dots, \SS_L)$ is realizable if and only if $\dim \SS_l \le r \cdot \dim \SS_{l-1}$ for $1 \le l \le L$. 

The initial lemma will be the basic tool for constructing module deformations with prescribed radical layerings. 

\proclaim{Lemma 5.1.  Constructing $\la$-module structures on graded vector spaces}  

\noindent We refer to the above notation.  Suppose that $\SS = (\SS_0, \dots, \SS_L)$ is a $d$-dimensional semisimple sequence.  Fix a direct sum decomposition of $K^d$ of the form 
$$K^d \ = \bigoplus_{0 \le l \le L} K_l, \ \ \  \text{where}\ \ \dim K_l = \dim \SS_l \ \ \ \text{for\ } l \le L,$$
and assume that $(f_\alpha)_{\alpha \in Q_1}$ is a family of $K$-linear maps $K^d \rightarrow K^d$, each taking $K_l$ to $\bigoplus_{u \ge l+1} K_u$.  {\rm{(}}Here $K_{L+1} = 0.${\rm{)}} Then $(f_\alpha)_{\alpha \in Q_1}$ is a point in $\modlasd$, and the corresponding $\la$-module $M$ has radical layering $\SS(M) \ge \SS$.  

Moreover, for any $m \in \{0, \dots, L\}$, the following conditions are equivalent:
\smallskip

$\bullet$ The first $m$ entries of $\,\SS(M)$ are $\, \SS_0, \dots, \SS_m$ in this order.
\smallskip
$\bullet$ For $l \in \{0, \dots, m\}$, the subspace $\sum_{q \in Q_{>0}} f_q (K_l)$ contains $\bigoplus_{u \ge l+1} K_u$.  
\medskip

In particular, $\SS(M) = \SS$ if and only if the following holds for $0 \le l < L$:  Modulo  $\bigoplus_{u \ge l+2} K_u$, the sum $\sum_{\alpha \in Q_1} f_\alpha(K_l)$ has the same dimension as $\SS_{l+1}$. 
\endproclaim 

\demo{Proof} The first assertion is immediate from the fact that $f_p(K^d) = 0$ whenever $p$ is a path of length $\ge L+1$.  To see that the $\la$-module $M$ corresponding to the family $(f_\alpha)_{\alpha \in Q_1}$ of linear maps satisfies $\SS(M) \ge \SS$, observe that, for any $l  \in \{1, \dots, L\}$, the $K$-vector space $J^l M$ equals 
$\sum f_p (K^d)$, where the sum extends over all paths $p$ of length $\ge l$; by hypothesis, this sum is contained in $\bigoplus_{j \ge l} K_j$.  Verification of the equivalences is straightforward.  \qed
\enddemo
  
The next lemma provides us with a means to exclude certain types of semisimple sequences from the set of potential generic radical layerings of irreducible components of $\modlasd$.

\proclaim{Lemma 5.2}  Let $\SS = (\SS_0, \dots, \SS_L)$ be a realizable semisimple sequence.  Suppose there exists an index $\rho \in \{0, \dots, L-1\}$ such that, for every $\la$-module $M$ in $\laySS$, the subfactor $J^\rho M/ J^{\rho +2} M$ has a simple direct summand.  Moreover, let $\SS_{\rho} = \SStilde_{\rho} \oplus  S$ and $\SStilde_{\rho + 1} = \SS_{\rho +1} \oplus S$, and suppose that the sequence
$$\SStilde = (\SS_0, \dots, \SS_{\rho - 1},\, \SStilde_{\rho},\,  \SStilde_{\rho + 1},  \SS_{\rho + 2}, \dots, \SS_L)$$ 
is again realizable.  Then $\laySS$ is contained in $\overline{\Mod \SStilde}$.
\endproclaim 

\demo{Proof} Let $M$ be any left $\la$-module with radical layering $\SS$.  Choose a basis 
$$B = (b_{l,j})_{0 \le l \le L,\, 1 \le j \le \dim \SS_l}$$
for $M$, with the property that, for $0 \le l \le L$, the elements $b_{l,1}, \dots, b_{l, \dim \SS_l}$ form a basis for $J^l M$ modulo $J^{l+1} M$.  The first hypothesis allows us to assume that $J\,b_{\rho,1}$ is contained in $J^{\rho + 2} M$.

In light of Observation 3.2, our realizability hypothesis guarantees that the semisimple sequence 
$(\SStilde_\rho, \SStilde_{\rho + 1}) = (S^{\dim \SS_{\rho} - 1}, \, \SS_{\rho + 1} \oplus S)$
is again realizable and that $\dim \SS_{\rho + 1} + 1 \le r \cdot (\dim \SS_{\rho} - 1)$.  Since $J^{\rho + 1} M / J^{\rho + 2} M = \SS_{\rho + 1}$ and $J b_{\rho, 1} = 0$ modulo $J^{\rho + 2} M$, this amounts to the existence of  an index $s \ge 2$ and an arrow $\gamma \in Q_1$ such that 
$$\gamma\, b_{\rho, s}\ \  \in \sum\Sb \alpha \in Q_1,\, j\ge2\\ \alpha \ne \gamma \text{\, or\, } j \ne s\endSb K\, \alpha\, b_{\rho, j} \ \ \ \text{modulo} \  J^{\rho + 2} M. \tag\dagger$$  
In particular, $\dim \SS_\rho \ge 2$.  It is clearly innocuous to assume $s = 2$.  

Using Lemma 5.1, we will construct a deformation $(D_t)_{t \in K}$ of $M$ in $\lamod$ such that $D_0 \cong M$, while $\SS(D_t) = \SStilde$ for all $t \ne 0$.  Our construction will amount to specifying a morphism $\AA^1 \rightarrow \modlasd$ which sends any $t \in K$ to a point $(f^{(t)}_\alpha)_{\alpha \in Q_1}$ in $\rep \SStilde$.  This will place $M$ into the closure of $\rep \SStilde$ in $\modlasd$, whence we will obtain the desired  containment $\laySS \subseteq \overline{\Mod \SStilde}$. 

We first pin down basis expansions of the products $\alpha \, b_{l,j}$ in $M$, where $\alpha$  traces the arrows of $Q$, namely
$$\alpha\, b_{l,j}\  = \sum_{u \ge l+1, \, v \le \dim \SS_u} c^{(\alpha,l,j)}_{u,v}\ b_{u,v}\ \ \ $$
with scalars $c^{(\alpha,l,j)}_{u,v}$.  For an unambiguous introduction of the $K$-linear maps $f^{(t)}_\alpha$, we rename the basis $B$ for $K^d$ to $\Btilde = (\btilde_{l,j})$.  In the spirit of Lemma 5.1 (relative to $\SStilde$), we introduce the vector space decomposition  
$$K^d \ = \bigoplus_{0 \le l \le L} K_l,$$
where the $K_l$ are as follows:  If $l \notin \{\rho, \rho+ 1\}$, we let $K_l$ be the subspace generated by the $\btilde_{l,j}$, $1 \le j \le \dim \SS_l$.  Moreover, $K_\rho$ is defined as the span of  $\btilde_{\rho,2}, \dots, \btilde_{\rho, \dim \SS_{\rho}}$, and $K_{\rho + 1}$ as the span of  the $(\dim \SS_{\rho + 1} + 1)$ elements $\btilde_{\rho, 1}$ and $\btilde_{\rho+ 1, j}$ for $1 \le j \le \dim \SS_{\rho + 1}$.

Our choice of $b_{\rho ,1}$ entails 
$$c^{(\alpha, \rho, 1)}_{\rho + 1, v} = 0\  \text{for all arrows}\   \alpha \ \text{and all}\  v \le \dim \SS_{\rho + 1}. \tag\dagger\dagger$$  
These equalities guarantee that the following $K$-endomorphisms $f^{(t)}_{\alpha}$ of $K^d$ satisfy the hypotheses of Lemma 5.1.  Namely, if $\alpha$ is an arrow in $Q_1 \setminus \{\gamma\}$, we define 
$$f_\alpha^{(t)} (\btilde_{l,j})\  =  \sum_{u \ge l+1, \, v \le \dim \SS_u} c^{(\alpha,l,j)}_{u,v}\ \btilde_{u,v}.$$ Moreover, we set 
$f_\gamma^{(t)} (\btilde_{l,j})  = \sum_{u \ge l+1, \, v \le \dim \SS_u} c^{(\gamma,l,j)}_{u,v}\ \btilde_{u,v}$ if either $l \ne \rho$ or $j  \ne 2$, and supplement by
$$f_\gamma^{(t)} (\btilde_{\rho, 2})\  =  \ \ t \cdot\btilde_{\rho, 1} \ \ + \ \ \sum_{u \ge \rho + 1, \, v \le \dim \SS_u} c^{(\gamma,\rho ,2)}_{u,\, v}\ \btilde_{u,v}.$$
Note that $f^{(t)}_\alpha(\btilde_{\rho, 1}) \in \bigoplus_{l \ge \rho + 2} K_l$ by $(\dagger\dagger)$, which shows that, for any arrow $\alpha \in Q_1$, we have  $f^{(t)}_{\alpha}(K_{\rho + 1}) \subseteq \bigoplus_{l \ge \rho + 2} K_l$. This makes the family of maps comply with the setup of Lemma 5.1.   Let $D_t$ be the $\la$-module defined by the maps $f^{(t)}_\alpha$.  That $D_0 \cong M$ is clear.   We now apply one implication of the final statement of Lemma 5.1 to the $K$-linear maps $K^d \rightarrow K^d, \ x \mapsto \alpha x$ defining $M$, and the other implication to the maps $f^{(t)}_\alpha$ defining $D_t$, in order to conclude that $D_t$ has radical layering $\SStilde$ whenever $t \ne 0$.  Indeed, from $(\dagger)$ and the construction of the $f_\alpha^{(t)}$ we deduce that 
$$\sum_{u \ge \rho + 1, \, v \le \dim \SS_u} c^{(\gamma,\rho ,2)}_{u,\, v}\ \btilde_{u,v} \ \ \in
\sum\Sb \alpha \in Q_1,\, j\ge2\\ \alpha \ne \gamma \text{\, or\, } j \ne 2\endSb K\, f^{(t)}_\alpha(\btilde_{\rho, j}) \ \ \ \ \text{modulo} \  \bigoplus_{l \ge \rho + 2} K_l\, ,$$
whence $\btilde_{\rho, 1}$ belongs to 
$$\sum_{j \ge 2}  f^{(t)}_{\gamma} (\btilde_{\rho,j})\ \  +  \ \sum_{\alpha \in Q_1,\, \alpha \ne \gamma,\, j \ge 2}  f^{(t)}_{\alpha} (\btilde_{\rho,j}) \ \ \ \ \text{modulo} \  \bigoplus_{l \ge \rho + 2} K_l$$ 
for $t \ne 0$.  Therefore $K_{\rho + 1}$ is contained in $\sum_{q \in Q_{>0}} f_q^{(t)} (K_{\rho})$.  We conclude that, for $t \ne 0$ and $l \ge 1$, the $K$-space $\sum_{\alpha \in Q_1} f^{(t)}_{\alpha}(K_{l - 1})$ equals 
$K_l$, modulo $\bigoplus_{u \ge l + 1} K_u$; but $\dim K_l = \dim \SStilde_l$ by definition.  As we pointed out earlier, this proves our claim. \qed
\enddemo

We remark that Lemma 5.2 does not apply to any $\rho \in \{0, \dots, L-1\}$ with the property that $S_{\rho}$ is simple, since the corresponding semisimple sequence $\SStilde$ fails to be realizable in this situation.  Hence the case of low dimensions $d$ will require separate consideration in the proof of Theorem A.

\demo{Proof of Theorem A}  We subdivide the argument into several steps.

{\bf{Step 1.}}  In this preliminary step, we let $d$ be any positive integer and suppose $\SS = (\SS_0, \dots, \SS_L)$ to be realizable  with $\dim \SS = d$.  We denote the closure of $\laySS$ in $\modlasd$ by $\U$ and set $\m = \min(d, L+1)$.  Our goal is to show that the irreducible variety $\U$  fails to be an irreducible component of $\modlasd$ unless $\SS_{\m - 1} \ne 0$.  

Assume $\SS_{\m - 1} = 0$, let $L'$ be maximal with respect to $\SS_{L'} \ne 0$, and pick $h \le L'$ maximal with respect to the condition that $\SS_h$ is not simple;  due to our assumption, such an index $h$ exists.  Given a module $M$ in $\laySS$, we choose a basis $B = (b_{l,j})$ for $M$, where $l \in \{0, \dots L'\}$ and $j \in \{1, \dots, \dim \SS_l\}$, such that the $b_{l,j}$ form a basis for $J^l M$ modulo $J^{l+1} M$.  For all arrows $\alpha$ and all legitimate choices of $l,j$, we then obtain expansions
$$\alpha\, b_{l,j} = \sum_{l+1 \le u \le L',\, v \le \dim \SS_u} c^{(\alpha,l,j)}_{u,v}\  b_{u,v}\ \  \tag\dagger$$
with $c^{(\alpha,l,j)}_{u,v} \in K$.  Since $\dim \SS_h \ge 2$, while $\dim \SS_{h+1} \in \{0, 1\}$ by our choice of $h$, we may assume  --  on possibly adjusting the basis elements $b_{h,j}$ for $J^hM$ modulo $J^{h+1}M$  --   that $\alpha_1 b_{h,1} \in J^{h+2} M$; indeed, since $\dim \SS_{h+1} = \dim J^{h+1} M / J^{h+2} M < \dim \SS_h$, the elements $\alpha_1 b_{h,j}$ with $1 \le j \le \dim \SS_h$ are linearly dependent modulo $J^{h+2}M$.   

Next we ascertain the existence of an index $s \ge 2$ with the property that $J^{h+1} M = J b_{h,s}$.  From Theorem 3.7 we glean that, generically, the modules $N$ in $\laySS$ (and hence the modules in $\U$) satisfy the following condition for each $l \le L$ with $\SS_{l+1} \ne 0$: For any arrow $\alpha \in Q_1$, there exists an element $x \in J^l N \setminus J^{l+1}N$ with the property that $\alpha\, x \in J^{l+1} N \setminus J^{l+2} N$.
In showing that $\laySS$ is properly contained in some closure $\overline{\Mod \SStilde}$ for a sequence $\SStilde <  \SS$ to be specified, it consequently suffices to verify that the modules $N \in \laySS$ with this latter property belong to the variety $\overline{\Mod \SStilde}$.  We return to the task at hand:  If $J^{h+1} M = 0$, there is nothing to be shown.  So suppose otherwise. $J^{h+1} M$ then being a nontrivial uniserial module, it suffices to ensure $J b_{h,s} \not\subseteq J^{h+2} M$ for some $s \ge 2$.  But in light of $\alpha_1 b_{h,1} \in J^{h+2} M$, the preceding considerations legitimize the assumption that $\alpha_1 \btilde_{h,s} \in J^{h+1} M \setminus J^{h+2} M$ for some $s \ge 2$.  It is clearly harmless to settle on $s = 2$.  In the sequel, we may thus assume 
$$J^{h+1} M = J b_{h,2}. \tag \dagger\dagger$$

Next we construct a family $(D_t)_{t \in K}$ of $d$-dimensional $\la$-modules such that $D_0  \cong M$ and $D_t$ has the following radical layering $\SStilde$ for $t \ne 0$:  Namely, $\SStilde_l = \SS_l$ for $l < h$, $\SStilde_h = S^{\dim \SS_h - 1}$, $\SStilde_{l} = S$ for $h+1 \le l \le L' +1$, and $\SStilde_l = 0$ for $l > L'+1$.  Note that, by our assumption, $L' + 1 \le L$.  

The sequence $\SS$ being realizable, Observation 3.2 shows the same to be true for $\SStilde$.  To unambiguously introduce the $D_t$, we rename the basis $B$ to $\Btilde = (\btilde_{l,j})$, cautioning, however, that the indexing of $\Btilde$ will not be in alignment with the radical layering $\SStilde$ of $D_t$ for $t \ne 0$;  the element  $b_{h,2}$ in the $h$-th radical layer of $M$ will be moved downward to turn into the reincarnation $\btilde_{h,2}$ in the $(h+1)$-st layer of $D_t$, and this will typically trigger further downward shifts.  In order to display a point $\bigl(f_{\alpha_i} \bigr)_{i \le r} \in \modlasd$ that gives rise to $D_t$ by way of Lemma 5.1, we consider the following direct sum decomposition of $d$-dimensional $K$-space: $K^d = \bigoplus_{0 \le l \le L'} K_l$, where $K_l =  \bigoplus_{j \le \dim \SS_l} K \btilde_{l,j}$ for $l < h$, $K_h = \bigoplus_{j \le \dim \SS_h,\, j \ne 2} K \btilde_{h,j}$, $K_{h+1} = K \btilde_{h,2}$, and $K_l = K \btilde_{l - 1,1}$ for $h+2 \le l \le L' +1$.  We start with the interloper assignment:
$$f^{(t)}_{\alpha_1} ( \btilde_{h,1})\ \  =\ \  t \cdot \btilde_{h,2} \ \ + \ \sum_{h + 2 \le u \le L'} c^{(\alpha_1, h, 1)}_{u,1}\  \btilde_{u,1}. \tag\dagger \dagger\dagger$$ 
To justify the range of the summation on the right of $(\dagger \dagger\dagger)$, keep in mind that $c^{(\alpha_1, h, 1)}_{h+1,v} = 0$ for $v \le \dim \SS_{h+1}$ by construction, and $\dim \SS_u = 1$ for $h+2 \le u \le L'$. 
Beyond $(\dagger\dagger\dagger)$, the linear maps $K_l  \rightarrow \bigoplus_{j \ge l+1} K_j$ formally duplicate those giving rise to the $\la$-module structure $M$ on $K^d$; see $(\dagger)$.  Explicitly, the $f_{\alpha}:  K^d \rightarrow K^d$ are defined via 
$$f^{(t)}_{\alpha}(\btilde_{l,j})\ \  =\ \   \sum_{u \ge l+1, \, v \le \dim \SS_u} c^{(\alpha,l,j)}_{u,v}\  \btilde_{u,v} \ \ \ \ \text{for} \ \ \ \ (\alpha,l,j) \ne (\alpha_1,h,1).$$
It is readily checked that the $r$-tuple $(f^{(t)}_{\alpha})_{\alpha \in Q_1}$ satisfies the conditions of Lemma 5.1, so as to yield a $\la$-module $D_t$ with $\SS(D_t) = \SStilde$ for $t \ne 0$.  Indeed, to verify this last equality, note that, for $t \ne 0$, equation $(\dagger \dagger \dagger)$ places $\btilde_{h,2}$ into $J^{h+1} D_t \setminus J^{h+2} D_t$;  then use $(\dagger\dagger)$ to show that the layers $\SS_l(D_t)$ for $l \ge h+2$ are as required.  It is obvious that $D_0 \cong M$. 

Clearly, the map $\AA^1 \rightarrow  \modlasd$ which sends $t$ to the point $\bigl(f_{\alpha_i}^{(t)} \bigr) \in \modlasd$ specified  above is a morphism of varieties.  This ensures that $M$ belongs to the closure of $\rep \SStilde$ in $\modlasd$, in turn an irreducible subvariety of $\modlasd$.  Given the $\U$-generic choice of $M$, we conclude that $\U \subsetneqq \overline{\Mod \SStilde}$.  Thus $\U$ fails to be an an irreducible component of $\modlasd$ as claimed.
\medskip

 {\bf{Step 2.}} In this step, we prove the (concluding) assertion of Theorem A,  addressing the case $d \le L+1$.  By Step 1, only the semisimple sequence $\SS$ with $\SS_l = S$ for $0 \le l \le d - 1$ and $\SS_l = 0$ for $l \ge d$ has the potential of being a generic radical layering of the modules in an irreducible component of 
 $\modlasd$.  Consequently, this sequence is the unique generic one for $\modlasd$, meaning that $\modlad = \overline{\laySS}$ is irreducible.  It is clear that the modules in $\laySS$ are precisely the uniserial modules of dimension $d$, whence the final statement of Theorem A is justified.  
 \medskip
 
{\bf{Step 3.}}  Finally, we verify the equivalence of conditions (1) -- (4) under the hypothesis that $d > L+1$.  The equivalence ``(1)$\iff$(1')" follows from Proposition 1.2.  Next we recall that realizability of $\SS$ is tantamount to $\dim \SS_l \le r \cdot \dim \SS_{l-1}$ for $1 \le l \le L$.

``(2)$\implies$(3)".  Assume that $(2)$ holds. Then $\SS$ is realizable by the first set of inequalities.  Let $E$ and $E_l$ be injective envelopes of $S$ and $\SS_l$, respectively.  Since $\soc_1(E) / S \cong S^r$, the companion inequalities, $r \cdot \dim \SS_l \ge  \dim \SS_{l-1}$ for $l \ge 1$, guarantee that the semisimple module $\SS_{l-1}$ embeds into the first socle layer $\SS^*_1(E_l) = \soc_1(E_l) / \soc_0 (E_l)$.  Using Lemma 2.8, we conclude that $(\SS_L, \dots, \SS_0)$ is indeed the generic socle layering of the variety $\laySS$. 

The  implication ``(3)$\implies$(4)" follows from  Lemma 2.5(a), and ``(4)$\implies$(1)" is covered by Theorem 3.1.

To verify ``(1)$\implies$(2)", we assume (2) to fail.  If $\dim \SS_l > r \cdot \dim \SS_{l-1}$ for some $l \in \{1, \dots, L\}$, then $\laySS$ is empty as remarked at the outset, and (1) fails as well.  So, in showing that the irreducible subvariety $\C = \overline{\laySS}$ of $\modlasd$ is not an irreducible component, we may assume $\dim S_l \le r \cdot \dim \SS_{l-1}$ for all eligible $l$.  Our assumption then yields an index $\rho \in \{0, \dots, L-1\}$ with the property that $\dim \SS_{\rho} > r \cdot \dim \SS_{\rho +1}\ \ (\dagger)$.  From Step 1, we moreover know that the equality $\SS_L = 0$ excludes the possibility that $\C = \overline{\laySS}$ be an irreducible component of $\modlasd$.  In ascertaining that $\C$ indeed fails to be such a component, it is therefore harmless to additionally assume $\SS_L \ne 0$.  

From $(\dagger)$ we infer that $\SS_{\rho}$ does not embed into $\soc_1(E_{\rho + 1}) / \soc_0(E_{\rho + 1})$, where $E_{\rho + 1}$ is an injective envelope of $\SS_{\rho + 1}$.  Consequently,  we obtain: $\bullet$ Every $\la$-module with radical layering $(\SS_{\rho}, \SS_{\rho+ 1})$ has a simple direct summand, and $\bullet$ the sequence 
$$\SStilde\  = \  (\SS_0, \dots, \SS_{\rho -1}, S^{\dim \SS_{\rho} - 1}, \SS_{\rho + 1} \oplus S, \SS_{\rho +2}, \dots, \SS_L)$$ 
is again realizable; keep in mind that $r \ge 2$.  Hence Lemma 5.2 implies that $\overline{\laySS}$ is (properly) contained in $\overline{\Mod \SStilde}$.  The latter variety being in turn irreducible, this rules out the possibility that $\overline{\laySS}$ is an irreducible component of $\modlasd$, and the argument is complete.
 \qed
\enddemo


\Refs
\widestnumber \key{ \bf 99}

\ref\no \BHT \by E. Babson, B. Huisgen-Zimmermann, and R. Thomas\paper Generic representation theory of quivers with relations \jour J. Algebra \vol 322 \yr 2009 \pages 1877--1918 \endref

\ref \no \BCH \by F. Bleher, T. Chinburg and B. Huisgen-Zimmermann \paper the geometry of algebras with vanishing radical square \paperinfo manuscript \endref 

\ref \no \BoHZtwo   \by K. Bongartz and B. Huisgen-Zimmermann \paper Varieties of uniserial representations  IV. Kinship to geometric quotients \jour Trans. Amer. Math. Soc. \vol 353 \yr 2001 \pages 2091--2113  \endref

\ref \no \CW \by A. T. Carroll and J. Weyman \paper  Semi-invariants for gentle algebras  \jour Contemp. Math. \vol 592 \yr 2013 \pages 111-136 \endref

\ref\no \CBS \by W. Crawley-Boevey and J. Schr\"oer \paper Irreducible
components of varieties of modules \jour J. reine angew. Math. \vol 553 \yr
2002 \pages 201--220  \endref

\ref\no\DoFl \by J. Donald and F. J. Flanigan \paper The geometry of Rep($A,V$) for a square-zero algebra \jour Notices Amer. Math. Soc. \vol 24 \yr 1977 \pages A-416 \endref

\ref\no\EiSa \by D. Eisenbud and D. Saltman \paper Rank varieties of matrices \inbook in Commutative Algebra (Berkeley 1987) \pages 173--212 \bookinfo MSRI Publ. 15 \publ Springer-Verlag \publaddr New York \yr 1989 \endref

\ref\no\Ger \by M. Gerstenhaber \paper On dominance and varieties of commuting matrices \jour Annals of Math. (2) \vol 73 \yr 1961 \pages 324--348 \endref

\ref\no\Gur \by R. M. Guralnick \paper A note on commuting pairs of matrices \jour Linear and Multilinear Algebra \vol 31 \yr 1992 \pages 71--75 \endref

\ref \no \hier \by B. Huisgen-Zimmermann  \paper A hierarchy of parametrizing varieties for representations \paperinfo in Rings, Modules and Representations (N.V. Dung, et al., eds.) \jour Contemp. Math. \vol 480   \yr 2009 \pages 207--239 \endref

\ref \no \GoHu \by B. Huisgen-Zimmermann and K. R. Goodearl \paper
Irreducible components of module varieties:  projective equations and rationality \jour Contemp. Math. \vol 562 \yr 2012 \pages 141--167 \endref

\ref \no \KacI \by V. Kac \paper Infinite root systems, representations of
graphs and invariant theory \jour Invent. Math. \vol 56 \yr 1980 \pages
57--92   \endref

\ref \no \KacII  \bysame \paper Infinite root systems, representations of
graphs and invariant theory, II \jour J. Algebra \vol 78 \yr 1982 \pages
141--162  \endref

\ref\no \Kra \by H.-P. Kraft \paper Geometric methods in
representation theory \inbook in Representations of Algebras (Puebla 1980) \eds M.
Auslander and E. Lluis \bookinfo Lecture Notes in Math. 944 \publ 
Spring\-er-Verlag \publaddr Berlin \yr 1982 \pages 180--258  \endref

\ref\no\Mor \by K. Morrison \paper The scheme of finite-dimensional representations of an algebra \jour Pac. J. Math. \vol 91 \yr 1980 \pages 199--218 \endref

\ref \no \RiRuSm \by C. Riedtmann, M. Rutscho, and S. O. Smal\o \paper Irreducible components of module varieties: An example \jour J. Algebra \vol 331 \yr 2011 \pages 130--144
 \endref
 
\ref \no \Scho \by A. Schofield \paper General representations of quivers
\jour Proc. London Math. Soc. (3) \vol 65 \yr 1992 \pages 46--64  \endref

\ref \no \Schro \by J. Schr\"oer  \paper Varieties of pairs of nilpotent matrices annihilating each other \jour Comment. Math. Helv. \vol 79 \yr 2004 \pages 396--426  \endref

\endRefs
\enddocument